\documentclass[12pt]{article}
\usepackage{amsmath,amssymb}%,sabon,hum521,euler}
\usepackage[T1]{fontenc}
\usepackage{graphicx}
\usepackage{amsthm}
\usepackage{manyfoot,perpage}
\DeclareNewFootnote{B}[fnsymbol]
\MakePerPage{footnoteB}

\theoremstyle{plain}
\newtheorem*{theorem}{Theorem}

\DeclareMathOperator{\disc}{disc}
\newcommand{\fp}{\mathfrak p}
\newcommand{\Of}{\mathcal O}
\newcommand{\Z}{\ensuremath{\mathbb{Z}}}
\newcommand{\Q}{\ensuremath{\mathbb{Q}}}
\newcommand{\F}{\ensuremath{\mathbb{F}}}
\newcommand{\C}{\ensuremath{\mathbb{C}}}

\let\phi\varphi
\DeclareMathAlphabet{\mathbb}{U}{msb}{m}{n}

\title{Dedekind on Higher Congruences and Index Divisors, 1871 and 1878.}
\author{Fernando Q.~Gouv\^ea and Jonathan Webster}
\date{}
\begin{document}

\maketitle

Dedekind's theorem connecting ideal theory and polynomial congruences
appears in all textbooks on algebraic number theory, but few books
note its connection to the problem of ``common index divisors.'' As
part of a project to study the history of this problem, we examine in
detail two of Dedekind's papers on the subject. (See also \cite{GW1}
for a similar analysis of Hensel's main work \cite{Hensel1894b} on the
same problem.)

We begin with a summary of the mathematical questions that Dedekind is
addressing here, then consider each of the publications in turn. In
each case we give a complete annotated translation. The first
publication is Dedekind's notice \cite{anzeige} of the second edition
of Dirichlet's \emph{Vorlesungen \"uber Zahlentheorie}. This naturally
focuses on Supplement X, the main addition to the text in the new
edition. Midway through, Dedekind announces some new results that were
proved and further elaborated in the second paper, \emph{\"Uber der
  Zusammenhang zwischen der Theorie der Ideale und der Theorie der
  h\"oheren Kongruenzen} \cite{Ded1878}. Both papers are translated
and annotated, both to clarify the mathematics and to highlight some
historical points.

% The first appears in a collection of
% publication notices by G\"ottingen faculty, the second in a kind of
% ``G\"ottingen proceedings'' that contains papers on various subjects
% by members of the University.

\section{Mathematical Background}

When Kronecker and Dedekind set out to generalize Kummer's theory of
cyclotomic integers, they quickly ran into obstacles. Finding a way
around these difficulties led each of them to develop a far more
complicated theory than Kummer's. As a result, each had to justify the
extra work by highlighting what made it necessary.

%% If we could be more precise here, we should.
Suppose $n>0$ is an integer and let $\zeta$ be a primitive $n$-th root
of unity. Kummer had found an explicit description in terms of
congruences of how rational primes factor in the cyclotomic integers
$\Z[\zeta]$. It seems that both Dedekind and Kronecker saw that
Kummer's description could be interpreted in terms of congruences
between polynomials (known as ``higher congruences'' at the time). In
modern terms, it would go something like this.

\begin{theorem} Let $n>2$ be an integer, let $\zeta$ be a
  primitive $n$-th root of unity, and let $\Phi_n(x)$ be the $n$-th
  cyclotomic polynomial. Fix a prime number $p\in\Z$ and let
  \[ \Phi_n(x)\equiv F_1(x)^{e_1} F_2(x)^{e_2} \dots F_r(x)^{e_r}
    \pmod p\] be the factorization of $\Phi(x)$ modulo $p$, where the
  $F_i(x)$ are distinct irreducible polynomials in $\F_p[x]$. Then the
  factorization of $(p)$ in $\Z[\zeta]$ is
  \[ (p)=\fp_1^{e_1}\fp_2^{e_2}\dots\fp_r^{e_r},\]
  with distinct prime ideals $\fp_i=(p,F_i(\zeta))$.
\end{theorem}

Of course, Kummer did not speak of ideals; instead, he thought of
$\fp_i$ as the ``ideal prime divisor'' determined by $p$ and
$F_i(x)$. This amounted to an explicit method for determining the
exponent of $\fp_i$ in a factorization. In modern terms, an ``ideal
prime divisor'' is essentially the valuation corresponding to $\fp_i$.

This beautiful result seemed to suggest the possibility of a very
simple theory in the general case: for a general number field
$\Q(\alpha)$, let $\Phi(x)$ be the minimal polynomial for $\alpha$ and
factor it modulo $p$. One could then use this to define ``ideal
primes'' \`a la Kummer.

The choice of $\alpha$ is crucial, of course. At least one example
would have been familiar to everyone: the field $\Q(\sqrt{-3})$ is the
same as the cyclotomic field of order $3$. Kummer's approach worked if
one took $\alpha$ to be a cube root of $1$ but would not work if we
took $\alpha=\sqrt{-3}$. Both Dedekind and Kronecker figured out that
one needed to work with all the algebraic integers in the field
$\Q(\alpha)$.
%% add refs? Eisenstein? Euler on FLT with n=3?

That highlights the first difficulty: in the case of $\Q(\zeta)$ the
ring of algebraic integers is exactly $\Z[\zeta]$, but this will not
be true in general. If $K$ is a number field and $\Of\subset K$ is its
ring of algebraic integers there may not exist any $\alpha\in\Of$ such
that $K=\Q(\alpha)$ and $\Of=\Z[\alpha]$. In such a situation, there
is no obvious $\Phi(x)$ to work with.

Under certain conditions we can still make it work, however. Given a
prime number $p\in\Z$, suppose we can find an $\alpha$ such that
$\Z[\alpha]\subset\Of$ has index not divisible by $p$. Then factoring
the minimal polynomial for $\alpha$ modulo $p$ gives the correct
factorization of $(p)$ in $\Of$. This theorem was announced by
Dedekind in \cite{anzeige} and proved in \cite{Ded1878}. It seems
clear that Kronecker was also aware of it.

This allowed one to hope, then, that an explicit factorization theory
could be based on a local approach: for each prime $p$, find a
generator $\alpha$ such that $p$ does not divide the index
$(\Of:\Z[\alpha])$. Then apply the theorem to find the
factorization. As we will see below, Dedekind says that he spent a
long time trying to prove that such an $\alpha$ always exists.

Alas, this is not true: there exist number fields in which \emph{all}
of the indices have a common prime divisor. Dedekind pointed this out
(and stated the factorization theorem) in \cite{anzeige}, probably to
explain why he had needed to take a different route. Kronecker says in
his \emph{Grundz\"uge} \cite[\S25, p.~384]{Grundzuge} of 1882 that he
had found an example in 1858.

%% Selling?
Both Dedekind and Kronecker pointed to this essential difficulty to
justify introducing a new approach: ideals in Dedekind's case, forms
in many variables in Kronecker's. Some years later, Zolotarev tried to
extend Kummer's theory directly in this style \cite{Zolo1}, but then
realized that his approach would fail for finitely many primes.
(Eventually, in a second paper \cite{Zolo2}, Zolotarev found still
another way to work around the difficulty.)  Dedekind's paper
\cite{Ded1878} was, as is clear from the introduction, prompted by
an announcement of Zolotarev's work.

Kronecker also mentioned Zolotarev's attempt in
\cite[\S25]{Grundzuge}, where he stated the problem in terms of
discriminants. For each choice of $\alpha$, let
$d(\alpha)=\disc(\Phi(x))$ be the discriminant of its minimal
polynomial. Let $d_K$ be the field discriminant. Then
$d(\alpha)=m^2 d_K$, where $m$ is exactly the index
$(\Of:\Z[\alpha])$. Kronecker, who always preferred specific elements
to collections, thought about this as follows. The many element
discriminants $d(\alpha)$ have a common factor $d_K$ which is the
essential part, attached to the ``Gattung'' $K$ rather than to a
specific element. The other factors of $d(\alpha)$ (i.e., the factors
of $m$) are therefore ``inessential.'' So in the ``bad'' examples what
is happening is that some prime $p$ is an inessential divisor of every
element discriminant. Such primes were the ``common inessential
discriminant divisors.''

The name is perhaps ill-chosen, because it is perfectly possible for a
prime $p$ to divide the discriminant $d_K$ and \emph{also} divide the
index $(\Of:\Z[\alpha])$. Such a prime divisor is then both
``essential'' (it divides $d_K$) and ``inessential''! Dedekind's
term ``index divisor'' seems more appropriate.

Kronecker's example ``in the thirteenth roots of 1'' is probably the
simplest one. He never gave the details, but they are probably as
Hensel gave them in his Ph.D.\ thesis \cite{HenselThesis} (see also
\cite[2.2]{Petri}). Let $\zeta$ be a primitive $13$-th root of
unity. There is a unique subfield $K$ of degree $4$ over
$\Q$.\footnote{This is global number field 4.0.2197.1 in
  \cite{LMFDB}.} Since the discriminant of $\Q(\zeta)$ is a power of
$13$, so is the discriminant of $K$ (in fact, $d_K=13^3$). It follows
from Kummer's work that the prime number $3$ is divisible by four
ideal primes in $K$, each of which has norm $3$; let $\fp$ be one of
these. Since $N(\fp)=3$, the field $\Of/\fp$ has three
elements.Consider some $\alpha\in K$. Since $K$ is a normal field, he
discriminant of the minimal polynomial of an integer in $K$ is the
product of differences of four integers in $K$. Since there are only
three congruence classes modulo $\fp$ at least one of these
differences must be divisible by $\fp$. Since $\fp$ lies above $3$,
the discriminant $d(\alpha)\in\Z$ must be divisible by $3$. Since
$d_K$ is a power of $13$, $3$ is an inessential divisor. In Dedekind's
terms, $3$ is a common inessential discriminant
divisor.

This set up the problem of determining exactly when this phenomenon
happens. One of the things that interests us about this problem is
that it was solved twice, apparently independently. Dedekind solved in
his paper \cite{Ded1878}. It is a sign of how little Kronecker
followed Dedekind's work that he suggested the problem of common
inessential discriminant divisors to Hensel for his Ph.D.\ in
1882. Hensel did not solve it completely in his
thesis,\footnote{Hensel later generalized the numerical condition in
  the example above to give a sufficient criterion for the existence
  of common inessential discriminant divisors, and even attempted to
  prove the condition was also necessary, which it is not. See the
  careful discussion in \cite[2.2]{Petri}.} but he published a
solution in 1894, in \cite{Hensel1894b}. While Hensel refers to
Dedekind's 1878 paper, it is not clear how carefully he had read
it. The relationship between the two solutions is complex and will be
explored elsewhere.

Two remarks might make it easier to read these texts. First, Dedekind
does not have the quotient construction, so when we would speak of
$\Of/\fp$ he must talk about congruence classes modulo $\fp$. For
example, the norm of an ideal is defined as the number of congruence
classes modulo that ideal. Second, despite the hints in Galois, at
this point there was essentially no theory of finite fields. Instead,
Dedekind relies on his paper \cite{Abriss}, in which he discusses
congruences modulo polynomials.

\section{Dedekind's \textit{Anzeige}}

The second edition of Dirichlet's \emph{Vorlesungen \"uber
  Zahlentheorie} \cite{DirDed2} appeared in 1871. The first edition,
which was in fact written by Dedekind based on his notes from
Dirichlet's lectures, had contained nine supplements, mostly taken
from Dirichlet's papers and supplementary lectures. (See
\cite{DirStill} for a translation (mostly) of the first edition; see
\cite{Landmarks} for more information on the book and its several
editions.)  The second edition included a new (tenth) supplement,
entitled ``On the Composition of Binary Quadratic Forms''
\cite{SuppX}. Most readers would have expected to find here a
simplified account of Gauss's theory. If they read it, however, they
would have been surprised to find, in the middle of the supplement, a
whole new theory of factorization in general fields of algebraic
numbers.

As he had done for the first edition in 1863, Dedekind wrote an
article \cite{anzeige} for the September 20, 1871 issue of the
\textit{G\"ottingische gelehrte Anzeigen}. This was a weekly journal
containing mostly book reviews, but Dedekind is of course writing
about his own book(s). As one would expect, the article focuses almost
entirely on the new content, i.e., Supplement~X. It gives a cursory
summary of the parts of the supplement containing well known material,
then discusses the new theory of ideals in more detail. Surprisingly,
after this explanation, Dedekind decided to go well beyond what is
found in Supplement X. In this digression, Dedekind announced results
about the relationship of ideal theory and ``higher congruences.''
These theorems were not in Supplement X and Dedekind would only
publish their proofs 1878 (\cite{Ded1878}, translated below). Once
this was done, he went back to a section-by-section discussion of the
Supplement.

Other authors had attempted to generalize Kummer's theory of ideal
factors in cyclotomic fields. Since Kummer had reduced his
factorization theory to congruence conditions, those attempts (like
the early attempts of Dedekind and Kronecker) had turned on ``higher
congruences.'' It seems that the discussion of ``common index
divisors'' (the name came later) was included in the book notice
exactly to explain why the approach via congruences was bound to
fail. This suggests that Dedekind was well aware that his
contemporaries would wonder whether something as innovative as the
theory of ideals was justified. His digression into ``higher
congruences'' makes the point that a new method is needed.

\subsection{Translation}

\textit{Vorlesungen \"uber Zahlentheorie}, by P.~G.~Lejeune
Dirichlet. Edited and with additions provided by R.~Dedekind. Revised
and enlarged second edition. Braunschweig, Friedr.\ Vieweg und
Sohn. 1871.

\vspace{\baselineskip}

I discussed the 1863 first edition of this work in these pages (27
January 1864), and so I can refer to that previous article with regard
to the origins and content of the book. The new edition differs from
the first by way of a great many additions, either in footnotes or to
the text itself. Many paragraphs\footnote{The book, including the
  supplements, is divided into 170 numbered ``paragraphs,'' which are
  actually sections, often several pages long. I will generally
  translate ``section'' from now on.} have been completely
reworked. These changes, which do not touch the essence of Dirichlet's
lectures, mainly reflect the decision to treat, in a new appendix, the
tenth supplement, the theory of the composition of binary quadratic
forms, which was omitted from the first edition for reasons discussed
at that time.\footnote{From here on everything in the book notice
  focuses on the content of Supplement~X.} The generality with which
Gauss presented the theory in the fifth section of the
\emph{Disquisitiones Arithmeticae} understandably causes significant
difficulties for the beginner. This difficulty motivated Dirichlet to
publish \textit{De formarum binariarum secundi gradus compositione},
1851.\footnote{This is \cite{DirFormarum}.} He says in his
introduction:

\begin{quote}
  De formarum compositione tunc non egi, quod argumentum ab
  illustrissimo Gauss in ``Disquisitionum Artihmeticarum'' sectione
  quinta maxima quidem generalitate sed per calculos tam prolixos
  tractatum esse constat, ut perpauci compositionis naturam percipere
  valuerint, eo magis quod summus geometra, ut ipse monuit, brevitati
  consulens theorematum difficuliorum demonstrationes synthetice
  adornavit, suppressa analysi per quam erant eruta. Quare confidere
  posse mihi videor, hujus argumenti expositionem novam et plane
  elementarem artis analyticae cultoribus non fore
  ingratam.\footnote{An idiomatic English translation might be: ``I
    did not then [in \cite{DirFormes}] take up the composition of
    forms, which topic was treated by the most illustrious Gauss in
    the fifth and largest section of this \textit{Disquisitiones
      Arithmeticae}. That treatment is so general and contains such
    long calculations that very few are able to grasp the nature of
    the composition. This is all the more so because that great
    geometer, as he himself admonished in keeping an eye for brevity,
    gave demonstrations of more difficult theorems by synthesis, with
    the analysis by which they were unearthed suppressed. For that
    reason it seems I may be confident that a new and undoubtedly
    elementary exposition of this theory will not be unwelcome to the
    cultivators of the analytic art.''}
\end{quote}

Since in this treatise only the first main theorem of the theory in
question is proved and no indication is given of how to continue, I
have taken a somewhat different route, which agrees with that of
Dirichlet in that only a special case of composition is
considered. Sections 145--149 contain the general theorems about the
composition of forms and classes of forms. These are used in sections
150, 151 to find the ratio of the class numbers of two determinants
whose ratio is a square; this is the same problem treated according to
Dirichlet's method in sections 97, 99, 100.\footnote{Dirichlet,
  following Gauss, assumes quadratic forms look like $ax^2+2bxy+cy^2$.
  Forms like $x^2+xy+y^2$ are replaced by $2x^2+2xy+2y^2$. As a
  result, they must allow forms where $\gcd(a,2b,c)\neq
  \gcd(a,b,c)$. This is the difficulty treated in sections 97--100.}
There follow in sections 152--154 the composition of genera and Gauss's
second proof of the quadratic reciprocity theorem. Sections 155--158
contain a proof of Gauss's theorem that the duplication of any class
results in the principal genus; it is based on a theorem of Lagrange
and Legendre about the rational integer solutions of indeterminate
equations of degree two in two unknowns.

In the paragraphs that follow I tried to introduce the reader to a
higher domain\footnote{This is the unexpected leap. To treat binary
  quadratic forms Dedekind wants to consider quadratic number
  fields. Since it is no harder (!) to treat all number fields, he
  proceeds to do so.} in which algebra and number theory are
intimately connected. In the course of my lectures on circle division
and higher algebra, held in G\"ottingen in winter 1856--1857 before
Mr.\ Sommer and Mr.\ Bachmann and in Winter 1857-1858 before Mr.\
Selling and Mr.\ Auwers, I was convinced that the study of the
algebraic properties of numbers is most appropriately based on
concepts that are directly linked to the simplest arithmetic
principles.\footnote{See \cite{Haubrich} for a reconstruction of the
  development of Dedekind's theory. See \cite{Haffner} for an analysis
  of what Dedekind means by ``simplest arithmetic principles.''}  I
replaced the name ``rational domain'' with the
word\footnote{Dedekind's work is of course \emph{K\"orper}, which
  translates as ``body.'' The standard English term is ``field,''
  which I will use throughout.}  ``field,'' by which I understand a
system of infinitely many numbers\footnote{By ``numbers'' Dedekind
  seems to mean complex numbers. So his fields are all subfields of
  $\mathbb C$.}  that has the property that the sums, differences,
products, and quotients of two such numbers belong to the same
system. I say a field $Q$ is a divisor of a field $M$, and the latter
a multiple of the first, when all the numbers contained in $A$ are
also found in $M$.\footnote{So $Q$ is a divisor of $M$ if
  $Q\subset M$.} Any two fields $A$, $B$ always have a least common
multiple, which can be denoted by $AB$, and also a greatest common
divisor.\footnote{We would call the lcm the compositum and the gcd the
  intersection of the two fields $A$ and $B$.} When to each number $a$
in a field $A$ there corresponds a number $b=\phi(a)$ so that
$\phi(a+a')=\phi(a)+\phi(a')$ and $\phi(aa')=\phi(a)\phi(a')$, the the
numbers $b$ make up a field $B=\phi(A)$ that is conjugate to $A$,
which arises from $A$ by the substitution\footnote{Dedekind uses
  ``substitution'' for what we would call a function, here a field
  homomorphism.} $\phi$. These concepts are connected, in the
algebraic direction, with the ideas of Galois and, in the
number-theoretical side, with Kummer's creation of the ideal
numbers.\footnote{As Dedekind said at the beginning of the paragraph,
  he is consciously creating a link between algebra and number theory,
  perhaps inspired by his teacher Dirichlet's linking analysis and
  number theory.}

In section 159 are developed the general properties of a field
$\Omega$ that has only a finite number of divisors.\footnote{Rather
  than use the dimension to define a finite extension, Dedekind
  focuses on the number of subfields, which allows him to stick to his
  ``sounds like arithmetic'' point of view. But he immediately points
  out that having finitely many subfields implies that there is a
  finite basis.} In such a field there is always a finite
quantity\footnote{We struggled to translate ``Anzhahl von Zahlen,''
  which is literally the awkward ``number of numbers.'' We settled for
  ``quantity of numbers.''} of numbers $\omega_1$, $\omega_2$,
$\dots$, $\omega_n$ with the property that any given number $\omega$
from the field can always and in only one way be expressed as
\[ h_1\omega_1+h_2\omega_2+\dots+h_n\omega_n,\] where $h_1$, $h_2$,
$\dots$, $h_n$ are rational numbers, which are called the coordinates
of the number $\omega$ with respect to the basis $\omega_1$,
$\omega_2$, $\dots$, $\omega_n$. The number $n$ is called the degree
of the field $\Omega$. It follows quite easily that every number in
the field is an algebraic number, namely the root of an equation of
degree $n$ whose coefficients are rational numbers. There are $n$
different substitutions linking the field $\Omega$ to conjugate
fields. The product of the $n$ values obtained from a given number
$\omega$ via these $n$ substitutions is called the norm of
$\omega$. It is a homogeneous function of the coordinates with
rational coefficients,\footnote{Dedekind remains attentive to the
  theory of forms, which is the ostensible subject of the supplement;
  the norm is a form of degree $n$ in $n$ variables.}  therefore a
rational number, which is denoted by $N(\omega)$. Given a system of
$n$ numbers $\alpha_1$, $\alpha_2$, \dots, $\alpha_n$ from the field
$\Omega$, one builds the determinant of the $n^2$ corresponding
numbers in the $n$ conjugate fields. The square of this determinant is
a rational number, which I call\footnote{This extension of the notion
  of discriminant seems to have appeared first here and in Supplement
  X.} the discriminant of the numbers $\alpha_1$, $\alpha_2$, $\dots$,
$\alpha_n$ and denote by $\Delta(\alpha_1, \alpha_2, \dots,
\alpha_n)$. It is not possible and not necessary to go into the
analytical developments\footnote{To a modern reader, the end of
  section 159 is very hard to follow; perhaps Dedekind's readers would
  have agreed. He skips all of it here.}  that are linked to these
concepts; they are only given in this paragraph to the extent that it
seemed appropriate for a better understanding.
%% look at what is in SuppX?

In the following section 160 all the algebraic numbers (which also
form a field) are divided into integral and fractional numbers. An
[algebraic] integer\footnote{Dedekind consistently writes ``ganze
  Zahl'' for an algebraic integer and ``ganze rational Zahl'' for an
  element of $\Z$. I will translate ``integer'' and ``rational
  integer'' respectively.} is understood to mean a root of an equation
with highest coefficient $=1$ and whose other coefficients are
rational integers. From this concept simple propositions
about divisibility, units, and relatively prime numbers are derived
for later use.

The following section 161 contains an auxiliary theorem for our theory
through which Gauss's notion of congruence between numbers can be
generalized.\footnote{This typical Dedekindian move feels perfectly
  comfortable for the modern reader, but it was not the way things
  were usually done in the 19th century. Dedekind here introduces a
  new algebraic idea, a ``module,'' and proceeds to establish the
  basic properties before returning to the theory of fields.} By a
module I understand a system $\mathfrak m$ of
numbers\footnote{Dedekind's modules are free $\Z$-submodules of $\C$.}
whose sums and differences still belong to the same system. The
congruence $\omega\equiv \omega'\pmod{\mathfrak m}$ means that the
difference $\omega-\omega'$ belongs to the system $\mathfrak m$ This
concept has a broader scope than its extraordinary simplicity seems to
promise, but we only give here what will serve to facilitate the
subsequent presentation.

After these preparations, the integers of a field $\Omega$ of degree
$n$ are investigated in section 162. They form a module $\Of$, and it
is shown first\footnote{The first step is to show the existence of an
  integral basis. A few lines later Dedekind will use the term
  ``fundamental series'' for such a basis.} that one can find $n$
integers $\omega_1$, $ \omega_2$, $ \dots$, $ \omega_n$ that are basis numbers
of the field, so that any integer can be represented as
\[ \omega=h_1\omega_1+h_2\omega_2+\dots+h_n\omega_n,\] where all the
coordinates $h_1$, $ h_2$, $ \dots$, $ h_n$ are whole numbers. The
discriminant $\Delta(\omega_1, \omega_2,\dots, \omega_n)$ of such a
basis, which I call a fundamental series,\footnote{Grundreihe.} has
the smallest possible absolute value. This nonzero rational integer is
of particular significance for the field $\Omega$, it is called the
discriminant or the fundamental number\footnote{Grundzahl.} and
denoted by $\Delta(\Omega)$. It divides the discriminant of any system
of $n$ integers, and the quotient is a square. Furthermore, if $\mu$
is a nonzero number in $\Of$, the number of incongruent integers with
respect to $\mu$ is equal to the absolute value of the norm
$N(\mu)$. We then draw attention to a strange phenomenon,\footnote{The
  ``strange phenomenon'' is the failure of unique factorization, but
  Dedekind describes it by saying that some irreducible elements of
  $\Of$ do not behave like true primes.} first observed in the case of
cyclotomic fields. It consists in this: a integer that cannot be
decomposed as a product of other integers does not always play the
role of a true prime number. This was the starting point for Kummer's
creation of ideal numbers.

My goal in section 163 is to propose a theory\footnote{A theory of
  factorization is meant. The ideal primes of Kummer will be replaced
  by prime ideals.} that applies to all [finite] fields. The
fundamental idea is as follows. If $\mu$ is a nonzero number in $\Of$,
then the system $\mathfrak m$ of all numbers in $\Of$ that are
divisible by $\mu$ has the following two properties:

I. The sum and difference of two numbers in $\mathfrak m$ is a number
in $\mathfrak m$; that is, $\mathfrak m$ is a module.

II. Every product of a number in $\mathfrak m$ with a number in $\Of$
is also a number in $\mathfrak m$.

It is not true that conversely, every system $\mathfrak m$ of integers
from a field that has these two properties, which from now on I will
call an ideal, is always the set of numbers that are divisible by some
fixed $\mu$. When this is the case, I say $\mathfrak m$ is a principal
ideal and denote it by the symbol $\mathfrak i(\mu)$.\footnote{In
  \cite{Ded1878}, the notation was changed to either $\Of\mu$ or
  $\Of(\mu)$, the latter when $\mu$ is an explicit number. See
  below.} We then investigate the properties of all the ideals of the
field $\Omega$, and the following main result follows. Multiplying
each number of an ideal $\mathfrak a$ by each number of an ideal
$\mathfrak b$, these products and their sums make up an ideal, which
is the product of the two factors $\mathfrak a$ and $\mathfrak b$ and
is denoted by $\mathfrak{ab}$.\footnote{Check against Supp X!} It then
clearly follows that $\mathfrak a\Of=\mathfrak a$,
$\mathfrak{ab}=\mathfrak{ba}$,
$(\mathfrak{ab})\mathfrak c = \mathfrak a (\mathfrak{bc})$, and that
from $\mathfrak{ab}=\mathfrak{ac}$ it follows that
$\mathfrak b=\mathfrak c$.  One says an ideal $\fp$ different from
$\Of$ is a prime ideal when it has no factors different from $\Of$ and
$\fp$;\footnote{Dedekind's definition of prime ideal sticks to the
  analogy with ordinary arithmetic. This definition is shown to be
  equivalent to the modern definition in \cite{Ded1878}; see
  p.~\pageref{decprime} below.} a composite ideal can be decomposed as
a product of prime ideals and in only one way. One then defines the
norm $N(\mathfrak a)$ of an ideal $\mathfrak a$ to be the quantity of
numbers in $\Of$ that are incongruent with respect to the module
$\mathfrak a$. We have
$N(\mathfrak{ab})=N(\mathfrak a)N(\mathfrak b)$. In this way we obtain
a complete analogy with the laws of divisibility in rational number
theory.

This entire theory is intimately connected with the so-called theory
of higher congruences,\footnote{Here Dedekind veers off the track. So
  far he has given a blow-by-blow account of Supplement X, but none of
  the material on higher congruences is found there.} which was
suggested by Gauss and developed om the work of Galois, Sch\"onemann
and others. It was first the works of Kummer on cyclotomic ideal
numbers and the study of the algebraic investigations of Galois that
led me to consider the theory of higher congruences, and I published a
brief outline of that theory (Crelle's Journal,
Vol.~54).\footnote{This is \cite{Abriss}, which discusses both
  congruences between polynomials and ``double congruences'' that
  amount to the theory of finite fields, hence the reference to
  Galois.} I later sought, with its help, to create a general theory
of ideal numbers, but was distracted from it by other work until the
publication of this\footnote{Presumably the publication of the first
  edition of the \textit{Vorlesungen}?} work led me back to that
subject. The renewed effort led me to my new theory of ideals, which
seems preferable to me because it is based on much simpler
concepts.\footnote{Many of Dedekind's contemporaries did not feel
  Dedekind's approach was in any way ``simpler.'' In particular, there
  was a lot of resistance to working with infinite sets as objects. As
  such, ideals seemed very abstract.}  In my presentation I did not
deal closely with the connection with the theory of higher
congruences, because I feared that the extent of my appendix would
become too large.\footnote{The supplement was 118 pages long, and,
  split into supplements X and XI, came to be much longer in later
  editions.} For readers who are interested in this
connection,\footnote{By which we suspect he means those who want to
  know why Dedekind did not stick to the straightforward approach.}  I
hope the following remarks may be useful.\footnote{Here begins the
  digression; this material is not in Supplement~X.}

Let $\omega$ be an arbitrary number in $\Of$, and
set\footnote{Dedekind does not mention that this is the discriminant
  of the polynomial of degree $n$ with $\omega$ as a root, but of
  course that was the original sense of ``discriminant'' that was
  generalized to $n$-tuples. He also does not mention the possibility
  that $D=0$; he is more explicit about this in the second paper.}
\[ \Delta(1,\omega,\omega^2,\dots,\omega^{n-1})=D^2\Delta(\Omega).\]
Then $D$ is always a rational integer, namely a homogeneous function
of degree $\frac12n(n-1)$ of the coordinates with rational integer
coefficients.\footnote{This homogeneous function was later called the
  ``index form.''} If then $p$ is a rational prime number and we are
given a number $\omega$ for which $D$ is not divisible by $p$, then
the decomposition of the principal ideal $\mathfrak i(p)$ as a product
of prime ideals is easily found via the theory of higher
congruences.\footnote{This is called ``Dedekind's Theorem'' in many
  modern textbooks. The proof was first given in \cite{Ded1878}; see
  below.} The number $\omega$ satisfies an equation of degree $n$
$F(\omega)=0$ and if
\[ F(x)\equiv P_1(x)^{e_1}P_2(x)^{e_2}\dotsc P_m(x)^{e_m} \pmod p,\]
where $P_1$, $ P_2$, $ \dots$, $ P_m$ are pairwise distinct prime
functions\footnote{Dedekind uses ``prime function'' for irreducible
  polynomial.} of the variable $x$ of degrees $f_1$, $f_2$, $\dots$,
$f_m$ respectively, then we have
\[ \mathfrak i(p) = \fp_1^{e_1}\fp_2^{e_2}\dotsc \fp_m^{e_m},\] where
$\fp_1$, $ \fp_2$, $ \dots$, $ \fp_m$ are pairwise distinct prime
ideals with norms $p^{f_1}$, $p^{f_2}$, $\dots p^{f_m}$,
respectively. From this follows easily\footnote{As Dedekind will
  clarify, it follows easily only for primes that do not divide $D$;
  see pages~\pageref{ramprimes1} and \pageref{ramprimes} below.} the
following theorem, which is fruitful for both algebraic and
number-theoretic investigations:

\textit{The prime number $p$ divides the fundamental number
  $\Delta(\Omega)$ of the field $\Omega$ if and only if $p$ is
  divisible by the square of a prime ideal.}\footnote{No italics in
  the original, but the statement does get its own paragraph.}

At first I thought it very likely\footnote{Compare the very similar
  comments in \cite{Ded1878} below, page~\pageref{longtime}.} that for
any given prime number $p$ there would exist an integer $\omega$ such
that the number $D$ was not divisible by $p$. Only when all my
attempts to prove the existence of such a number were unfruitful did I
set myself the task of investigating whether this conjecture was
incorrect.\footnote{As indeed it is, which Dedekind will show.
  Kronecker claimed that he knew this in 1858.} Were the conjecture
true, whenever $p$ is divisible by $r$ distinct prime ideals $\fp$
whose norms have value $p^f$, there must exist $r$ distinct prime
functions $P$ of degree $f$. Conversely,\footnote{This is one of the
  main results in \cite{Ded1878}; it is proved again in
  \cite{Hensel1894b}.} when this last condition is always satisfied,
then one can prove the existence of a number $\omega$ with the desired
property. In the simplest case when $f=1$, there are exactly $p$
distinct prime functions of degree one. The question then becomes
whether there exists a field $\Omega$ in which $p$ is divisible by
$(p+1)$ distinct prime ideals, all of which of norm $p$. The degree of
such a field must then be $=p+1$. The simplest case arises when one
takes $p=2$, leading to the question: do there exist cubic fields in
which the number $2$ is divisible by three distinct prime ideals? In
such a field $D$ would always be an even number.\footnote{Dedekind has
  twice reduced to the ``simplest case'' in order to find his
  example. The task now is to find a cubic field in which $2$ splits
  completely, which will force $D$ to be even for every choice of
  $\omega$. In his example, however, Dedekind proves \emph{direclty}
  that $D$ is always even by computing it explicitly, and then appears
  to \emph{conclude} that $2$ splits completely.} One can always
assume that the fundamental series\footnote{Dedekind will begin by
  taking an integral basis and considering the corresponding
  multiplication table. This provides him with a number of parameters
  he can adjust to obtain the desired field.} of a cubic field
consists of the number $1$ and two integers $\alpha$, $\beta$ whose
product is rational.\footnote{If $\alpha\beta=\ell\alpha+m\beta+n$,
  replace $\alpha$ by $\alpha-m$ and $\beta$ by $\beta-\ell$.} One
then has\footnote{Dedekind offers no explanation for why the formulas
  should look like this. To spare the reader some time, here is an
  explanation. Define integers $a$, $a'$, $b$, $b'$, $c$, $c'$, $n$ by
  $\alpha\beta=n$, $\alpha^2=a'\alpha+b\beta-c$ and
  $\beta^2=a\alpha+b'\beta-c'$. Then notice that
  $n\beta=\alpha\beta^2$. Expanding the latter and equating basis
  coefficients gives $n=ab$, $c=bb'$, $c'=aa'$, as Dedekind
  says. Note, however, that $a$, $b$, $a'$, $b'$ are not arbitrary:
  the minimal polynomials of $\alpha$ and $\beta$ depend on them, and
  bad choices will give polynomials that are not irreducible, so that
  the $\Q$-algebra defined by these equations will not be a field.}
\begin{align*}
  \alpha\alpha &= a'\alpha+b\beta-bb'\\
  \beta\beta&= a\alpha + b'\beta-aa'\\
  \alpha\beta &= ab
\end{align*}
where $a$, $b$, $a'$, $b'$ are rational integers with no common
divisor,\footnote{If some prime divides all four integers, then $p^2$
  would divide $\alpha^2$, $\alpha\beta$, and $\beta^2$, and hence
  $p^2$ would divide $(\alpha+\beta)^2$, and so
  $\frac1p\alpha+\frac1p\beta\in\Of$, contradicting the assumption
  that $\{1,\alpha,\beta\}$ is an integral basis.} and we can
compute\footnote{Given the information we have, we can compute the
  traces of $\alpha^2$, $\beta^2$, and $\alpha\beta$; from that
  information it is easy to compute the discriminant.}
\[ \Delta(\Omega)=\Delta(1,\alpha,\beta)\]
\[ = {a'}^2{b'}^2+18aba'b'-4a{a'}^3-4b{b'}^3-27a^2b^2.\]
If we now set
\[ \omega= z+x\alpha +y\beta,\]
with $z,x,y$ any rational integers, then
\begin{multline*}
  \omega^2=z^2 = z^2 - bb'x^2-aa'y^2-aa'y^2+2abxy\\
  +(a'x^2+ay^2+2xz)\alpha + (bx^2+b'y^2+2yz)\beta,
\end{multline*}
and it follows\footnote{We have expressed the basis
  $\{1,\omega,\omega^2\}$ as a linear combination of
  $\{1,\alpha,\beta\}$; $D$ is the determinant of that matrix.} that
\[ D = bx^3-a'x^2y+b'xy^2-ay^3,\] independent of $z$, which is
expected from the definition of $D$. Even though $a$, $b$, $a'$, $b'$
have no common divisor, $D$ will be an even number whenever $a$ and
$b$ are even and $a'$ and $b'$ are odd.\footnote{If $a,b$ are even and
  $a',b'$ are odd,
  \[ D\equiv x^2y+xy^2\equiv xy(x+y)\equiv 0 \pmod{2}\] for all
  integers $x$, $y$. Thus, Dedekind has shown directly that $2$ is a
  common index divisor for any cubic field of this form. It remains to
  show that there is actually a choice of $a$, $b$, $a'$, $b'$ that
  makes it all work.} It must then be that the number $2$ is divisible
by three distinct prime ideals. This is completely confirmed by the
example\footnote{Dedekind now chooses the quadruple $(2,2,1,-1)$. To see
  that this is not a random choice, note that the minimal polynomial
  for $\alpha$ is $x^3 - a'x^2 + bb'x - ab^2$. While this is
  irreducible for most choices of the quadruple $(a,b,a',b')$, that is
  not always the case. Dedekind's choice gives
  $x^3-x^2-2x-8$, which is indeed irreducible and so we have a cubic
  field of discriminant $-503$, global number field 3.1.503.1 in
  \cite{LMFDB}. But if we chose $(2,2,1,1)$ we would get
  $x^3-x^2+2x-8=(x-2)(x^2+x+4)$. Even more dramatically, $(6,2,9,13)$
  would give $\Delta=1$, which is impossible for a number field, and
  indeed $x^3-9x^2+26x-24 =(x-2)(x-3)(x-4)$.}
\[ a=b=2,\quad a'=-b'=1,\quad \Delta(\omega)=-503;\] we
have\footnote{To confirm that he has the example he wants, Dedekind
  writes out (without proof) the factorizations of
  $2,\alpha,\beta$. Note, however, that he has already shown that $2$
  is a common index divisor; given that, his converse theorem forces
  the factorization to be as he wants, since there do exist
  irreducible polynomials of degree two and three in $\F_2[x]$.}
\[ \mathfrak{i}(2)=\mathfrak{abc},\quad
  \mathfrak{i}(\alpha)=\mathfrak{a^2c}, \quad
  \mathfrak{i}(\beta)=\mathfrak{b^2c},\] where
$\mathfrak{a,\, b,\, c}$ are three distinct prime ideals.\footnote{See
  the more complete discussion in \cite{Ded1878},
  page~\pageref{details} below, where Dedekind defines the three
  ideals explicitly and checks all of these statements.}

Another example can be obtained in the following way. With respect to
the modulus $p=2$ there exists only one prime function of degree two,
namely $x^2+x+1$. Therefore when in a field $\Omega$ the integer $2$
is divisible by at least two distinct prime ideals whose norm
$=p^2=4$, then $D$ must be even. In this case the degree of the field
must be at least $=4$. The phenomenon in fact occurs in the
biquadratic field\footnote{This is global number field 4.0.2873.1 in
  \cite{LMFDB}.} defined by the equation
\[ \alpha^4-\alpha^3+\alpha^2-2\alpha+4=0.\]
The numbers $1$, $\alpha$, $\beta=2:\alpha$, and
$\gamma=\alpha^2-\alpha$ are a fundamental series and the fundamental
number is $=13^2\cdot17$.

Thus, there exist fields $\Omega$ in which the number $D$ above is
always divisible by certain singular prime numbers $p$. Of course
there are only finitely many such primes. I remark, however, that the
theorem\label{ramprimes1} above, characterizing of the rational primes
that divide the fundamental number $\Delta(\Omega)$ of a field,
remains valid in general, but it would take us too far afield were I
to give a proof of this theorem or to explore its significance for the
theory of fields.

After this digression, I continue to summarize the contents of the
sections that follow. In section 164 the ideals of the field $\Omega$
are divided into a finite number of classes. Two ideals are called
equivalent when their product by some fixed ideal is a principal
ideal. An ideal class consists of all ideals that are equivalent to a
given ideal. The principal class consists of the principal
ideals. These ideal classes allow a composition that has the same
properties of the composition of classes of quadratic forms.

In section 165 I show the relationship between the composition of
ideal classes and the decomposable homogeneous forms that arise from
the same field $\Omega$.\footnote{These are the forms given by the
  norm function; they are ``decomposable'' because the norm is a
  product by definition.}

Section 166 gives Dirichlet's theory of units in a generalized
form. The presentation is completely independent of the previous
one. In section 167 this theory is used to obtain an expression for
the number of ideal classes by way of an infinite series, much like
the determination of the class number of quadratic forms. At this
point I break away from the study of the general problem, since my
investigations of this topic have not yet been crowned with
sufficient success to be published. The sections that follow,
168--170, illustrate the general theory by applying it to the example
of quadratic fields.

So far it appears that the theory of ideal numbers has been the
subject of serious research by only four or five
mathematicians.\footnote{Who were they?} My heartfelt wish is that the
new edition of Dirichlet's \textit{Vorlesungen \"uber Zahlentheorie}
may facilitate access to this large subject and perhaps to motivate a
larger number\footnote{Famously, this did not happen. See, for
  example, the correspondence with Lipschitz translated in
  \cite[Section 0.7]{B}.} of mathematics to apply their powers so
that, amidst the tremendous recent progress made in geometry and in
the theory of functions, number theory may not be left behind.

\vspace{\baselineskip}

\hspace*{0.5in} July 22 1871 \hspace*{2in} R.~Dedekind

\section{The 1878 Paper}

The title of \cite{Ded1878} translates as ``On the Relationship
between the Theory of Ideals and the Theory of Higher Congruences.''
It tells us that the paper will discuss the connections between two
subjects about which Dedekind had already written: the theory of
``higher congruences'' and the theory of ideals. By ``higher
congruences'' Dedekind means not only congruences modulo a prime
between polynomials of higher degree but also the kind of congruence
he will write as ``$\mathrm{modd}~p,P$,'' where $p$ is a prime and
$P$ is a polynomial.\footnote{Today we would describe this as working
  in the quotient ring $\F_p[x]/(P)$.} The notation
``$\mathrm{modd}$'' is intended to call attention that there are two
moduli in play.\footnote{We will nevertheless write
  $\mathrm{mod}~p,P$.}  Higher congruences had been discussed by
Dedekind in his \emph{Abri\ss} of 1857, which includes much of what we
would now describe as the theory of finite fields. That paper is one
of the main references used; Dedekind denotes it as ``C'' for short.

At this point, Dedekind had given two accounts of his theory of
ideals: first in \cite{SuppX}, Supplement X of the second edition of
Dirichlet's \emph{Vorlesungen \"uber Zahlentheorie} and then in an
article \cite{B} published in French. Dedekind refers to these
as ``D'' and ``B''; for an English reader B is the preferred
reference, since it was translated by John Stillwell and published by
Cambridge University Press \cite{B}.

In what follows we give a loose annotated translation of
\cite{Ded1878}. The translation is ``loose'' in the sense that we have
not tried to preserve the exact syntactic structure of Dedekind's long
sentences nor always attempted (and certainly not always succeeded) to
capture every nuance of meaning. We have, however, tried to translate
the mathematical content precisely, mostly preserving Dedekind's
terminology. Our annotations are given as numbered footnotes;
Dedekind's own footnotes are marked with asterisks. Page numbers in
\cite{Gesam} are indicated in the margin. Dedekind numbers his
main results as I, II, III, etc.; we have labeled those theorems
accordingly, but have highlighted other results (usually given by
Dedekind in italics) as theorems as well.

In \cite{Gesam}, \"Oystein Ore added several endnotes, which we give
in summary form at the end. The editors also added a few footnotes
that we have translated in annotations, distinguishing them from
Dedekind's original footnotes. There are several spelling changes made
in \cite{Gesam}; for example, ``Discriminante'' becomes
``Diskriminante.'' When we quote the German, we have tried to stick to
the original spelling.

\subsection{Translation}

\begin{center}
  {\large On the Relationship between the Theory of
    Ideals\\[.2\baselineskip] 
    and the Theory
  of Higher Congruences}\\[\baselineskip]
by R.~Dedekind
\end{center}

\vspace{2\baselineskip}

\marginpar{[202]} The new principles by which I arrived at a theory of
ideals that is rigorous and without exceptions were first explained
seven years ago in the second edition of the \emph{Lectures on Number
  Theory} by Dirichlet (\S~159--170) and more recently given, in
greater detail and in slightly modified form, in the \emph{Bulletin
  des sciences math\'ematiques et astronomiques} (t.~XI. p.~ 278; t.~I
(2e. serie), p.~17, 69, 144, 207).\footnote{The first publication is
  \cite{SuppX}, 1871, and the second is \cite{B}, 1876. The
  third edition of Dirichlet's \emph{Vorlesungen}, which contained a
  version of the theory similar to \cite{B}, was yet to appear.}
Stimulated by the great discovery of Kummer, I had been concerned with
this subject for many years, starting from a completely different
basis, namely the theory of higher congruences.  Although these
investigations brought me very close to the desired goal, I decided
not to publish them, because the theory that emerges suffers from two
imperfections.  The first is that the investigation of a domain of
integral algebraic numbers begins first with the consideration of a
certain number and the equation corresponding to it, which is then
interpreted as a congruence. The definitions of ideal numbers (or
rather of divisibility by ideal numbers) are obtained in this
way. Since everything depends on a specific representation, it follows
that the invariant character of the definition cannot be recognized
from the start.\footnote{This is one of Dedekind's fundamental
  methodological principles: one should always try to define things in
  a way that is independent of specific choices, rather than making
  such choices and then proving invariance. He wanted his mathematics
  ``coordinate-free.''} The second imperfection of this approach is
that there are peculiar exceptional cases that require special
treatment.\footnote{Those special cases are the index divisors, one of
  the main topics of this paper.} My more recent theory, on the other
hand, is based exclusively on notions such as \emph{fields},
\marginpar{203]} [algebraic] \emph{integers}, and \emph{ideals}, whose
definition does not require any particular form of representation of
the numbers, removing the first defect. The power of these extremely
simple concepts is shown by the fact that, in proving the general laws
of divisibility, a distinction between several cases never occurs
again.  I have made some remarks about the connection between the two
types of justification and stated some theorems without proof in the
\textit{G\"ottingischen gelehrten Anzeigen} of September 20, 1871
(pp. 1488--1492). In particular I have discovered the
reason\footnote{Perhaps Hensel wrote \cite{Hensel1894b} because he did
  not think Dedekind's ``reason'' was a sufficient answer to the
  question; see \cite{GW1}.}  for the existence of the
peculiar exceptional cases mentioned above. Since then, a theory of
ideal numbers by Zolotareff appeared in 1874, in a paper in Russian
with the title \textit{Th\'eories des nombres entiers complexes, avec
  une application au calcul integral}.\footnote{This is \cite{Zolo1};
  the ``completion'' mentioned below eventually appeared
  as \cite{Zolo2}.}  This was announced and abstracted in the
\textit{Jahrbuch \"uber die Fortschritte der Mathematik} (Vol. 6, p
117). From the abstract\footnoteB{I can only refer to the
  abstract. After several unsuccessful attempts to get it in the
  bookstore, I have recently obtained the original through the
  kindness of Professor Wangerin, but given my ignorance of the
  Russian language, to my great regret I was able to understand very
  little, only what is clear from looking at the formulas.} it is
clear that the theory of Zolotareff is also based on the theory of
higher congruences, but that the treatment of the aforementioned
exceptional cases is temporarily excluded and is reserved for a later
presentation.  I do not know if this prospective completion has since
been published.  Since, however, the connection between the two types
of justification of general ideal theory is of sufficient interest in
itself, I allow myself to provide here the proofs of the remarks given
in the \textit{G\"ottingischen gelehrten Anzeigen}.\footnote{It seems,
  then, that Zolotarev's paper was the main stimulus for writing this
  paper.}

I will assume as known both my theory of ideals and the theory of
higher congruences, of which I gave a short description earlier in
Borchardt's \textit{Journal} (Vol.~54, p.~1).\footnote{This is
  \cite{Abriss}; Borchardt was then the editor of the \textit{Journal
    f\"ur die Reine und Angewandte Mathematik}.} For brevity, I will
cite this paper on congruences as C, the second edition of Dirichlet's
number theory\footnote{This is \cite{DirDed2}, but more specifically
  \cite{SuppX}.}  as D, and the paper in the \textit{Bulletin des
  sciences math\'ematiques}\footnote{This is \cite{AlgInts}, but we
  cite the English translation \cite{B}.} as B.

\vspace{2\baselineskip}
\centerline{\S~1\footnote{This section introduces the key
  objects in play: the ring of integers $\Of$ of a number field
  $\Omega$, the order $\Of'=\Z[\theta]$, and the index $k$.}}
\vspace{2\baselineskip}

\marginpar{[204]} Let $\Omega$ be a finite field\footnote{Dedekind
  says ``finite field'' for what we would call ``a finite extension of
  \Q.'' He never considers fields with finitely many elements.} of
degree $n$, and let ${\Of}$\footnote{Dedekind uses the lowercase
  fraktur $\mathfrak o$.} be the domain of all [algebraic]
integers\footnote{Dedekind uses ``ganzen Zahlen,'' literally ``whole
  numbers,'' for algebraic integers. I will typically translate
  ``integers.'' The elements of $\Z$ are ``rational integers.''}
contained in it. There always exist $n$ independent
integers \[\omega_1,\omega_2,\dots,\omega_n\] which are a basis for
the domain $\Of$, that is, the system $\Of$ is identical with the
collection
\[ [\omega_1,\omega_2,\dots,\omega_n]\]
of all numbers $\omega$ of the form
\[ \omega=h_1\omega_1+h_2\omega_2+\dots+h_n\omega_n.\]
where
\[ h_1,h_2,\dots,h_n\] are arbitrary rational integers.  The
discriminant\footnote{The notion of the discriminant of a set of
  algebraic numbers seems to have been created by Dedekind by analogy
  to the older notion of the discriminant (or determinant) of a
  polynomial. In this paper Dedekind typically uses ``Discriminante''
  for the general construct, reserving ``Grundzahl,'' which we
  translate as ``fundamental number,'' for this particular
  discriminant.}
\[ \Delta(\omega_1,\omega_2,\dots,\omega_n)=\Delta(\Omega)=D,\] which
is independent of the choice of the basis numbers\footnote{Dedekind
  uses ``Basiszahlen.'' From here on I will use the modern term
  ``integral basis.''}  $\omega_1,\omega_2,\dots,\omega_n$, is called
the fundamental number or the discriminant of the field
$\Omega$. (D.~\S~159, 160,162; B.~\S~12--18).

Now if $\theta$ is a specific algebraic integer in the field, we can
set\footnote{Dedekind writes $c_i'$ where I have $c_i^{(1)}$, $c_i''$
  where I have $c_i^{(2)}$, etc.}
\begin{align*}
  1 &= c_1^{(0)}\omega_1 + c_2^{(0)}\omega_2+ \dots +c_n^{(0)}\omega_n\\
  \theta &= c_1^{(1)}\omega_1 + c_2^{(1)}\omega_2+ \dots +c_n^{(1)}\omega_n\\
  \theta^{2} &= c_1^{(2)}\omega_1 + c_2^{(2)}\omega_2+ \dots +c_n^{(2)}\omega_n\\  
  \dots &= \dots\\
  \theta^{n-1} &= c_1^{(n-1)}\omega_1 + c_2^{(n-1)}\omega_2+ \dots
  +c_n^{(n-1)}\omega_n
\end{align*}
where all the $n^2$ coefficients or coordinates $c$ are rational
integers, and we will have 
\[\Delta(1,\theta,\theta^2,\dots,\theta^{n-1})=Dk^2,\]
where\footnote{The meaning of this notation, standard at the time, is
  $k=\det[c_i^{(j)}]$. Dedekind does not use matrices, which had not
  yet been invented, nor does he represent the determinant as an
  array.}
\[ k=\sum\pm c_1^{(0)}c_2^{(1)}\dots c_n^{(n-1)}\] is a rational
integer. The absolute value of this number $k$, which is independent
of the choice of integral basis, will for brevity from now one be
called the \emph{index} of the integer $\theta$. If $k$ is not $0$, as
we will always assume,\footnote{This running assumption is crucial,
  but it is not mentioned every again. It is equivalent to assuming
  that $\theta$ is a generator.} the numbers
\[1,\theta,\theta^2,\dots,\theta^{n-1}\] will be 
independent of each other (D.~\S~159; B.~\S~4,15,17) and $\theta$ will
be the root of an irreducible  equation of degree $n$
\[ F(\theta)=\theta^n+a_1\theta^{n-1}+\dots + a_n=0,\] where the
coefficients $1,a_1,a_2,\dots,a_n$ are all rational integers.

\marginpar{[205]}If we let $\phi(t)$ be any function of the variable
$t$, --- and I remark that always, by this name [function] and by an
expression of the form $\phi(t)$, $f(t), \dots$ in this treatise one
should always understand an entire function of $t$ whose coefficients
are rational integers\footnote{So ``function'' always means a
  polynomial with integer coefficients.} --- the set $\Of'$ of all
other numbers of the form \[\omega'=\phi(\theta)\] is called an
\emph{order} (D.~\S~165, 166; B.~\S~23); all such numbers are integers
of the field $\Omega$ and therefore are contained also in
$\Of$. Clearly it suffices to take only the
functions \[\phi(t)=x_0+x_1t+x_2t^2+\dots+x_{n-1}t^{n-1}\] whose
degree is smaller than $n$, since if $\phi_1(t)$ has degree larger
than $n$ we can divide it by
\[F(t)=t^n+a_1t^{n-1}+a_2t^{n-2}+\dots+a_{n-1}t+a_n.\] The remainder
$\phi(t)$ will have degree less than $n$ and at the same time
$\phi_1(\theta)=\phi(\theta)$.  In the notation used above (B.~\S~3)
we can set
\[ \Of'=[1,\theta,\theta^2,\dots,\theta^{n-1}].\] It also follows from
the irreducibility of the equation $F(\theta) = 0$ that each number
$\omega'$ can be represented in the form
$\phi(\theta)$\footnote{Dedekind means $\phi(t)$ with $\phi$
of degree less than $n$.} in only one way; nevertheless, in what
follows we will not always restrict ourselves to that form of
representation, but rather allow functions of any degree.

Prime numbers $p$ --- by which name we mean a rational positive prime
number --- fall in two cases once the fixed number $\theta$ is
chosen:\footnote{The idea, then, is to choose and fix $\theta\in\Of$
  such that $\Omega=\Q(\theta)$. Then $\Of'=\Z[\theta]\subset \Of$ and
  $k=(\Of:\Of')$ is the index. The rational primes $p$ that divide $k$
  are those in the second case; the (infinitely many) other primes are
  in the first case.}  the \emph{first} case, which applies to
infinitely many prime numbers, is when the index $k$ of the number
$\theta$ is \emph{not} divisible by $p$.  If $k=\pm 1$, then all
primes are in this first case, and in fact $\Of'$ is identical to
$\Of$. When however $k^2>1$, a finite number of primes will fall into
the \emph{second} case, namely the prime divisors of $k$. The
paragraphs that follow will show that the decomposition of the prime
numbers $p$ of the first kind\marginpar{[206]} (or rather the
decomposition of the corresponding principal ideals\footnoteB{This
  notation for principal ideals is more appropriate than
  $\mathfrak{i}(p)$, which I used earlier (D.~\S~163).}  $\Of p$) as a
product of prime ideals\footnote{Here and elsewhere Dedekind writes
  ``Produkte aus lauter Primidealen,'' literally ``product of nothing
  but prime ideals,'' to emphasize that it is a complete factorization
  into primes.}  can be completely reduced to the decomposition of the
function $F(t)$ as a product of functions that are prime\footnote{We
  would say irreducible modulo $p$.} with respect to the modulus $p$
(C.~6).\footnote{This is the theorem announced in the \emph{Anzeige}
  \cite{anzeige}.} On the other hand, it is not possible to do this in
the same simple way for prime numbers of the second kind. The
following remarks should be made before this investigation.

Let $p$ be a fixed prime of the \emph{first} kind, so that $k$ is not
divisible by $p$. In this case an element of
$\Of'$
\[ \omega' = x_0+x_1\theta+ x_2\theta^2 + \dots +
  x_{n-1}\theta^{n-1}\] is divisible by $p$ (which means it is equal
to $p\omega$, with $\omega$ an integer, i.e., an element of $\Of$) if
and only if all the coefficients $x_i$ are divisible by $p$. To see
that, notice\footnote{This paragraph is simple linear
  algebra. Dedekind uses the expression of the powers of $\theta$ in
  terms of the integral basis to rewrite $\omega'$ in terms of the
  basis.} that
\[ \omega'=h_1\omega_1+h_2\omega_2+\dots+h_n\omega_n\]
where
\begin{align*}
  h_1 &=c_1^{(0)} x_0+c_1^{(1)} c_1+c_1^{(2)} x_2+\dots +c_1^{(n-1)} x_{n-1}\\
  h_2 &=c_2^{(0)} x_0+c_2^{(1)} c_1+c_2^{(2)} x_2+\dots +c_2^{(n-1)} x_{n-1}\\
\dots & \dots\\
  h_n &=c_n^{(0)} x_0+c_n^{(1)} c_1+c_n^{(2)} x_2+\dots +c_n^{(n-1)} x_{n-1}
\end{align*}
It follows from the independence of $\omega_1,\omega_2,\dots,\omega_n$
that $\omega'$ is divisible by $p$ if and only if each of the
coordinates $h_1,h_2,\dots,h_n$ is divisible by $p$. If so, each of
the products $kx_0,kx_1,kx_2,\dots,kx_{n-1}$ is also divisible by $p$,
and therefore so are the coefficients
$x_0,x_1,x_2,\dots,x_{n-1}$.\footnote{If the matrix $[c_i^j]$ is
  invertible mod $p$, then $\vec h\equiv\vec 0$ if and only if
  $\vec x\equiv\vec 0$.}  The same theorem can clearly also be stated
as: a number $\omega'$ of the order $\Of'$ is divisible by a prime
number $p$ of the first kind if the quotient $\frac{\omega'}{p}$ is
itself in the order $\Of'$. Conversely, when all the coefficients
$x_0,x_1,x_2,\dots,x_{n-1}$ are all divisible by $p$, then obviously
$\omega'$ is divisible by $p$.\footnote{Added a paragraph break here.}

Therefore\footnote{Dedekind will now translate congruences between
  elements $\phi(\theta)$ of $\Z[\theta]$ into congruences between the
  polynomials $\phi(t)$. He does so first under the assumption that
  $\deg(\phi)<n$, and then in the general case. For the latter he uses
  the ``double modulus'' $p,F(t)$.} two numbers $\phi_1(\theta)$ and
$\phi_2(\theta)$ of the order $\Of'$ are congruent modulo $p$ (i.e.,
their difference $\phi_1(\theta)-\phi_2(\theta)$ is divisible by $p$)
if and only if the coefficients of the two functions $\phi_1(t)$ and
$\phi_2(t)$ are all congruent modulo $p$, i.e., in the sense of the
theory of higher congruences, when we have
\[ \phi_1(t)\equiv \phi_2(t)\pmod p\]
(C.~1). For this conclusion, however, we need to assume that the
degrees of the functions $\phi_1(t)$ and $\phi_2(t)$ are less than
$n$. If that is not the case, after dividing by $F(t)$ we obtain an
identity of the form
\[\phi_1(t)-\phi_2(t)=F(t)\psi(t)+\psi_1(t),\]
where $\psi_1(t)$ has degree less than $n$, and then
$\phi_1(\theta)-\phi_2(\theta)=\psi_1(\theta)$. We will have
\[ \phi_1(\theta)\equiv \phi_2(\theta) \pmod p\] when
$\psi_1(t)=p\psi_2(t)$, that is, when
\[\phi_1(t)-\phi_2(t)=F(t)\psi(t)+p\psi_2(t).\] The existence of such
an identity is described in the theory of higher congruences as
\[\phi_1(t)-\phi_2(t)\equiv F(t)\psi(t)\pmod p \] or simply as (C.~7)
as \[\phi_1(t)\equiv \phi_2(t) \pmod{p,F(t)}.\] Conversely, it is
clear that from that function congruence the number
congruence \[\phi_1(\theta)\equiv \phi_2(\theta)\] also follows; the
two congruences are therefore equivalent. Thus in $\Of'$ there are as
many numbers $\phi(\theta)$ that are incongruent modulo $p$ as there
are functions $\phi(t)$ incongruent with respect to the double modulus
$p, F(t)$; there are $p^n$ of the latter (C.~8), which is also the
number\footnote{Dedekind writes $(\Of,\Of p)$ for the index.}
$(\Of:\Of p)=N(p)$ of numbers in $\Of$ that are incongruent modulo $p$
(B.~\S~18; D.~\S~162), which implies the following result:
\textit{each number $\omega$ of the domain $\Of$ is congruent modulo
  $p$ to a number $\omega'$ of the order $\Of'$.}\footnote{What
  Dedekind has shown is that when $p$ does not divide the index the
  quotient $\Of/p\Of$ is isomorphic to $\Of'/p\Of'$. As he shows next,
  if $k$ is the index $(\Of:\Of')$ the isomorphism is given by
  multiplication by $lk$, where $l$ is any rational integer such that
  $lk\equiv 1\pmod p$.}

The same conclusion can be reached directly by the following simple
argument. From the $n$ relations between the numbers
$1,\theta,\theta^2,\dots,\theta^{n-1}$, on the one hand, and the
numbers $\omega_1,\omega_2,\dots,\omega_n$, on the other, it follows
that the products $k\omega_1,k\omega_2,\dots,k\omega_n$ are contained
in the order $\Of'$, and therefore so are all the products $k\omega$
for any $\omega$ belonging to $\Of$. Therefore we have
$k\omega=\phi(\theta)$. Now since $k$ is not divisible by $p$, we can
choose a rational integer $l$ such that $kl\equiv 1\pmod p$, and then
it follows that $\omega\equiv lk\omega\equiv l\phi(\theta) \pmod p$,
so that $\omega$ is congruent modulo $p$ to a number $l\phi(\theta)$
that is in the order $\Of'$.

Things are completely different when $p$ is a prime
the\marginpar{[208]} \emph{second} kind.\footnote{For index divisors,
  $\Of'/p\Of'\hookrightarrow \Of/p\Of$ is not onto.} In that case the
determinant $k$ is divisible by $p$, and it is easy to prove that
there exist $n$ rational integers $x_0,x_1,\dots, x_{n-1}$, not all
divisible by $p$, such that the corresponding numbers
$h_1,h_2,\dots, h_n$ are all divisible by $p$. Then the corresponding
number
\[ \omega'=x_0+x_1\theta + x_2\theta^2 + \dots +x_{n-1}\theta^{n-1}\]
is in fact divisible by $p$ even though the coefficients
$x_0,x_1,\dots, x_{n-1}$ are not all divisible by $p$. It follows that
the number $(\Of':\Of p)$ of incongruent elements in $\Of'$ is smaller
than $p^n$.\footnote{We would write $(\Of':\Of'p)$, but Dedekind does
  not.} It follows that there are numbers $\omega$ in $\Of$ that
are not congruent modulo $p$ to \emph{any} element $\phi(\theta)$,
i.e., there exist congruence classes $(\text{mod}~p)$ in $\Of$ for
which there is no representative in $\Of'$. The precise determination
of the number $(\Of':\Of p)$ is not necessary for our
purposes.\footnote{The editors of \cite{Gesam} add a footnote here:
  ``In Zolotareff one also finds the theorem that the exceptional
  prime numbers are precisely those for which there is a number
  $\omega'$ in the order $\Of'$ that are divisible by $p$ but whose
  coefficients are not all divisible by $p$. Zolotareff does not say,
  however, that these prime numbers are the index divisors.''}

% Editors' footnote: Schon bei Zolotareff findet man den Satz, dass dei
% Ausnahmeprimzahlen eben diejenigen sind, wofur eine durch p teilbare
% Zalh omega' in der Ordnung O' vorkommt, worin nicht alle Koeffiziented
% durch p telbar sind. Zolotareff seight abeg nich, dass diese
% Primzahlen eben die Indexteiler sind.

\vspace{2\baselineskip}
\centerline{\S~2}
\vspace{2\baselineskip}

%\subsection{Paragraph 2: the Main Theorem}

In this paragraph\footnote{This section states and proves
  ``Dedekind's Thoerem,'' describing the factorization of primes of
  the first kind in terms of higher congruences. Dedekind will use the
  phrase ``prime function'' to mean a monic polynomial that is
  irreducible modulo $p$. Everything in this section assumes that
  $p$ does not divide $k$.} we consistently make the assumption that
$p$ is a prime number of the \emph{first} kind. We want to prove
that in this case the theory of higher congruences gives an easy way
to decompose a principal ideal $\Of p$ into its prime factors. This
happens because the function $F(t)$, which we will denote $F$ for
brevity, factors modulo $p$ as a product of \emph{prime functions}
$P(t)$ (C.~6). If we assume, for convenience, that each prime function
$P$ has highest coefficient $=1$, it follows that two incongruent
prime functions are always relatively prime (C.~5). Combining all the
congruent factors into powers we get
\[ F \equiv P_1^{e_1} P_2^{e_2} \dots P_m^{e_m} \pmod p\] where the
$P_i$ are all the incongruent prime functions contained in
$F$.\footnote{Dedekind knows that there is unique factorization in
  $\F_p[x]$. This is one of the many results in \cite{Abriss}.}

Let $P$ be any one of these $m$ prime functions, and \marginpar{[209]}
let $\rho=P(\theta)$. Then there is an ideal $\mathfrak p$ that is the
greatest common divisor of $\Of p$ and $\Of \rho$. To study the
properties of this ideal $\mathfrak p$, we first determine all the
elements $\phi(\theta)$ contained in the order $\Of'$ that are
divisible by $\mathfrak p$ (i.e., are contained in $\fp$). We want to
prove\footnote{Dedekind doesn't state this as a separate Lemma but he
  uses it over and over in the sequel.}  that the congruence
\[ \psi(\theta) \equiv 0 \pmod \fp \tag{1}\]
is completely equivalent to the function congruence
\[ \label{psit} \psi(t)\equiv 0 \pmod{p,P} \tag{2}.\]
Indeed,\footnote{Here begins the proof. Recall that $\rho=P(\theta)$
  where $P$ is an irreducible factor of $F$.} by definition (D.~\S~163;
B.~\S~19) the ideal $\fp$ is the collection of all numbers of the form
\[ \rho\alpha+p\beta,\]
where $\alpha,\beta$ are arbitrary numbers from the domain $\Of$. By
\S~1, each number $\alpha$ is congruent modulo $p$ to some number
$\phi(\theta)$ in the order $\Of'$, so from (1) we get a congruence of
the form
\[ \psi(\theta) \equiv P(\theta)\phi(\theta) \pmod p;\]
  this is equivalent (as in \S~1) to the function congruence
\[\psi(t)\equiv P(t)\phi(t)\pmod{p, F},\]  
and therefore also equivalent to congruence (2), since $F$ is
divisible by $P$.  Conversely, it follows immediately\footnote{From
  (2) we get $\psi(t)=P(t)\phi(t) + pG(t)$; plugging in $\theta$ gives
  $\psi(\theta)=\rho\phi(\theta)+pG(\theta)\in\fp$, since
  $\phi(\theta)\equiv\alpha\pmod p$.}  from (2) that
any $\psi(\theta)$ is of the form $\rho\alpha+p\beta$, and so is
$\equiv 0 \pmod\fp$ as well. This proves our claim above.\footnote{The
  lemma is now proved.}

With the help of these results we can easily\footnote{The argument is
  to pass from $\Of/\fp$ to $\Of'/\fp$ and then to translate
  congruences between elements of $\Of'$ into congruences of
  polynomials using the Lemma above. This reduces the problem to
  counting incongruent polynomials modulo $p,P$, which Dedekind had
  already done in \cite{Abriss}.} compute the \emph{norm}
of the ideal $\fp$, i.e., the number $(\Of:\fp)=N(\fp)$ of elements of
$\Of$ that are incongruent modulo $\fp$. So let $\alpha_1,\alpha_2$ be
any two numbers in $\Of$. From \S~1 we now that there exist two numbers
$\phi_1(\theta)$, $\phi_2(\theta)$ in $\Of'$ that are congruent modulo
$p$ to $\alpha_1,\alpha_2$. Since $\fp$ divides $p$, we also have
\[ \alpha_1\equiv\phi_1(\theta),\quad \alpha_2\equiv\phi_2(\theta)
  \pmod\fp.\] So the two numbers $\alpha_1,\alpha_2$ are congruent
modulo $\fp$ if and ony if
\[ \phi_1(\theta)\equiv \phi_2(\theta) \pmod\fp.\]
This congruence is equivalent, as above, to the congruence
\[ \phi_1(t)\equiv \phi_2(t) \pmod{p,P}.\] Therefore there are as many
numbers $\alpha$ that are incongruent modulo $\fp$ as there are
functions $\phi(t)$ incongruent with respect to the
\marginpar{[210]}double modulus $p,P$; this quantity is $=p^f$, where
$f$ is the degree of the function $P$ (C.~8), so we have
\[ N(\fp)=p^f.\]

With that it is easy to prove $\mathfrak p$ is a \emph{prime
  ideal}. First, we know $f\geq 1$, so $N(\mathfrak p)\neq 1$, so
$\mathfrak p$ cannot be equal to $\Of$. It suffices then to
show\footnote{Here we see that (at this time) Dedekind's working
  definition of ``prime ideal'' is not the same as the one we learn
  today. The next several lines explain why it is enough to prove that
  $\alpha_1\alpha_2\in \mathfrak p$ implies that either $\alpha_1$ or
  $\alpha_2$ is in $\mathfrak p$. The argument is straightforward;
  note that Dedekind consistently writes ``$\mathfrak m$ divides
  $\alpha$'' instead of ``$\alpha$ belongs to
  $\mathfrak m$.''\label{decprime}} that $\fp$ is not a decomposable
ideal, i.e., that it is not a product of the form
$\mathfrak a_1\mathfrak a_2$, where $\mathfrak a_1\mathfrak a_2$ are
ideals and neither is equal to $\Of$. Such a decomposable\footnote{The
  proof starts here.} ideal $\mathfrak m= \mathfrak a_1\mathfrak a_2$
has the characteristic property that there are always two numbers
$\alpha_1,\alpha_2$, neither divisible by $\mathfrak m$, whose product
$\alpha_1\alpha_2$ is divisible by $\mathfrak m$. This is because both
the ideals $\mathfrak a_1,\mathfrak a_2$ are different from $\Of$, so
neither of them can be divisible by their product
$\mathfrak m= \mathfrak a_1\mathfrak a_2$. So there must exist a
number $\alpha_1$ that is divisible by $\mathfrak a_1$ but not by
$\mathfrak m$, and similarly an $\alpha_2$ that is divisible by
$\mathfrak a_2$ but not by $\mathfrak m$. So $\fp$ will be a prime
ideal if we can show\footnote{We have shown that an ideal $I$ is
  indecomposable if and only if $ab\in I$ implies either $a\in I$ or
  $b\in I$.} that a product $\alpha_1\alpha_2$ cannot be divisible by
$\fp$ unless at least one of the factors $\alpha_1,\alpha_2$ is
divisible by $\fp$. For this,\footnote{Now we will prove $\fp$ is
  prime; as usual, we reduce to elements of $\Of'$ and then to
  polynomial congruences.} we set, as above,
\[ \alpha_1\equiv \phi_1(\theta),\alpha_2\equiv
  \phi_2(\theta)\pmod\fp,\] 
so that
\[ \alpha_1\alpha_2\equiv \phi_1(\theta)\phi_2(\theta)\pmod\fp,\] and
since $\alpha_1\alpha_2\equiv 0\pmod\fp$, we must have  
\[ \phi_1(\theta)\phi_2(\theta)\equiv 0 \pmod{\mathfrak p}\]
and so
\[ \phi_1(t)\phi_2(t)\equiv 0 \pmod{p,P}.\] Since $P$ is a \emph{prime
  function} it follows\footnote{Irreducibles in $\F_p[t]$ are prime,
  which Dedekind had proved in \cite{Abriss}.}
that one of the two congruences
\[ \phi_1(t)\equiv 0 \pmod{p,P} \qquad \text{ or }
\qquad\phi_2(t)\equiv 0 \pmod{p,P}\]
must hold (C.~6). So at least one of the congruences
\[ \phi_1(\theta)\equiv 0 \pmod{\fp} \qquad \text{ or } \qquad
\phi_2(\theta)\equiv 0 \pmod{\fp}\] must be true, that is, one of the
two numbers $\alpha_1,\alpha_2$ must $\equiv 0\pmod\fp$. Therefore 
$\mathfrak p$ is a prime ideal, and we know (B.~\S~21) that $\fp$ is a
prime ideal of \emph{degree} $f$, since $N(\mathfrak p)=p^f$.

Now\footnote{We have a prime ideal dividing $p$, so it remains to
  determine the valuation, i.e., the highest power of $\fp$ dividing
  $p$.} we would like to prove that the highest exponent $e$ of $P$ in
the factorization of $F$ and the highest exponent of $\mathfrak p$ in
the factorization of $p$ are equal. Indeed, if $F$ is
divisible modulo $p$ by $P^e$ but not by $P^{e+1}$, We have
\[F \equiv S P^e\pmod p,\] where $S$ is not divisble by
$P$. \marginpar{[211]} It follows as above that \[\sigma=S(\theta)\]
is not divisible by $\mathfrak p$. Since $\mathfrak p$ is the greatest
common divisor of $\Of p$ and $\Of \rho$, we know that
\[ \Of p = \mathfrak{pa},\qquad \Of \rho=\mathfrak{pb}\] with
$\mathfrak a$ and $\mathfrak b$ relatively prime. So what we need to
prove\footnote{Dedekind first proves that $\mathfrak a$ is divisible
  by $\fp^{e-1}$, and then proves it cannot be divisible by
  $\fp^{e}$.} is that the highest power
of $\mathfrak p$ contained in $\mathfrak a$ is $\mathfrak p^{e-1}$.
For this, consider the number
\[ \eta=\sigma\rho^{e-1} = S(\theta)P(\theta)^{e-1},\] which cannot be
divisible by $p$, since the degree of the polynomial $SP^{e-1}$ is
less than $n$ and its highest coefficient is $=1$. On the other hand,
$\eta$ is divisible by $\fp^{e-1}$, since $\rho$ is divisible by
$\mathfrak p$. From the congruence $F\equiv SP^e\pmod p$, we see that
$\eta\rho=\sigma\rho^e$ is divisible by $p$. So the ideal
$\eta\mathfrak{pb}$ is divisible by $\mathfrak{pa}$, and therefore
$\eta\mathfrak b$ is divisible by $\mathfrak a$; since $\mathfrak a$
and $\mathfrak b$ are relatively prime, we see that $\eta$ is
divisible by $\mathfrak a$. So let \[\Of \eta=\mathfrak{ac},\] where
$\mathfrak c$ is an ideal not divisible by $\fp$,\footnoteB{It follows
  that $\mathfrak a$ is the greatest common divisor of the ideals
  $\Of p$ and $\Of \eta$, and so $\eta\fp$ is the least common
  multiple of $\Of p$ and $\Of \eta$, i.e., $\fp$ is the collection of
  all roots $\pi$ of the congruence $\eta\pi\equiv 0 \pmod p$. This
  could also have been used to define the ideal $\fp$.} because
otherwise $\eta$ would be divisible by $\mathfrak{ap}=\Of p$, which we
know is not the case. Since $\eta$ is divisible by
$\mathfrak p^{e-1}$, so is $\mathfrak a$.\footnote{Added a paragraph
  break here. The first part of the proof is finished: $\mathfrak a$
  is divisible by $\fp^{e-1}$ and so $\Of p=\fp\mathfrak a$ is
  divisible by $\fp^e$; in valuation terms,
  $v_\fp(\mathfrak a)\geq e-1$.}

We now only need to show that $\mathfrak a$ is not divisible by
$\mathfrak p^e$.  Since\footnote{The argument opens with ``if
  $\mathfrak a$ is divisible by $\fp^e$,'' which seems like setting up
  a proof by contradiction. But that is not where the argument goes. I
  think it is better expressed as two cases: if $\mathfrak a$ is not
  divisible by $\fp$, then since $0=v_\fp(\mathfrak a)\geq e-1$ we
  must have $e=1$ and we are done. If $\mathfrak a$ is divisible by
  $\fp$, then Dedekind shows that $v_\fp(\rho)=1$ and so the number
  $\eta=\sigma\rho^{e-1}$ is divisible by $\mathfrak a$ but not by
  $\fp^e$. From that it follows that $v_\fp(\mathfrak a)\leq e-1$,
  hence must be exactly $e-1$.}  $e\geq 1$, if $\mathfrak a$ is
divisible by $\mathfrak p^e$, then it is certainly divisible by
$\mathfrak p$ itself. Now if $\mathfrak a$ is divisible by $\fp$,
$\mathfrak b$ cannot be divisible by $\fp$, and therefore $\rho$ is
not divisible by $\mathfrak p^2$. From that it follows that $\sigma$
is not divisible by $\mathfrak p$, so in this case $\mathfrak p^{e-1}$
is the highest power of $\mathfrak p$ contained in the number
$\eta=\sigma\rho^{e-1}$. So $\eta$, and therefore the ideal
$\mathfrak a$ contained in it, cannot be divisible by $\mathfrak p^e$,
which was to be proved.

After this the investigation of a specific\footnote{For each different
  irreducible factor $P$ we have found a prime ideal $\fp$ dividing
  $p$ and shown that the multiplicity of $P$ as a factor of $F$ is the
  same as the multiplicity of $\fp$ as a factor of $p$. To complete
  the proof we put these all together and then show that these are no
  other prime ideals dividing $p$.} prime function $P$ contained in $F$
and its corresponding prime ideal $\fp$ is complete. We now apply the
results to all the functions contained in $F$,
\[ F\equiv P_1^{e_1}P_2^{e_2}\dots P_m^{e_m} \pmod p,\]
with incongruent prime functions \marginpar{[212]}
\[ P_1, P_2, \dots, P_m\]
of degrees, respectively,
\[ f_1, f_2,\dots, f_m.\]
To these functions correspond prime ideals
\[ \mathfrak p_1,\mathfrak p_2, \dots, \mathfrak p_m\]
with the corresponding degrees, so that
\[ N(\mathfrak p_1)=p^{f_1}, N(\mathfrak p_2)=p^{f_2}, \dots,
  N(\mathfrak p_m)=p^{f_m}\] 
and
\[ \mathfrak p_1^{e_1},\mathfrak p_2^{e_2},\dots, \mathfrak
  p_m^{e_m}\] are the highest powers of these ideals contained $p$.
These $m$ ideals are all distinct;\footnote{This is a key fact later
  on: distinct factors correspond to distinct ideals and vice
  versa. Dedekind shows only that different $P$ give different
  ideals.}  for example, $P_2$ is not divisible by $P_1$ mod $p$, so
the number $P_2(\theta)$ is divisible by $\mathfrak p_2$ but not by
$\mathfrak p_1$, and it follows that $\mathfrak p_1$ and
$\mathfrak p_2$ are different ideals.  Finally,\footnote{The last
  thing to note is that we have the complete factorization: no other
  prime ideals divide $p$. This is easy to see.} we know that $p$
cannot be divisible by any other prime ideal, since\footnote{The
  product of all the $P_i(\theta)^{e_i}$ is congruent mod $p$ to
  $F(\theta)=0$.}
\[P_1(\theta)^{e_1}P_2(\theta)^{e_2}\dots P_m(\theta)^{e_m}\equiv
  0 \pmod p.\] If $p$ is divisible by a prime ideal, that
ideal has to divide one of the $m$ numbers $\rho=P(\theta)$; but then
that ideal must be identical to the prime ideal $\mathfrak p$, which
is the greatest common divisor of $\Of p$ and $\Of\rho$.

From all this it follows (D.~\S~163, B.~\S~25) that
\[\Of p = \mathfrak p_1^{e_1}\mathfrak p_2^{e_2}\dotsc
  \mathfrak p_m^{e_m}.\] A consequence\footnote{Surprisingly, this
  famous formula does not appear in \cite{SuppX} or \cite{B}.}
of this, found by taking norms, is
\[ n = e_1f_1+e_2f_2+\dots +e_mf_m.\] Thus we have proved the
following theorem,\footnote{Usually known today as ``Dedekind's
  theorem.''} which I announced in the \emph{G\"ottingischen gelehrten
  Anziegen} in September 20, 1871.

\begin{theorem}[I] Let $k$ be the index\footnote{Recall
    the running assumption that $k\neq 0$, so that $\Omega=\Q(\theta)$
    and the miminal polynomial is of degree $n$.}
  of the number $\theta$ that satisfies the irreducible equation of
  degree $N$ $F(\theta)=0$. If $k$ is not divisible by $p$
  and if
  \[ F \equiv P_1^{e_1}P_2^{e_2}\dots P_m^{e_m} \pmod p\] where the
  $P_1,P_2,\dots,P_m$ are incongruent prime functions of degree
  $f_1,f_2,\dots, f_m$, respectively, then we have
  \[\Of p = \mathfrak p_1^{e_1}\mathfrak p_2^{e_2}\dotsc
    \mathfrak p_m^{e_m},\] where
  $\mathfrak p_1,\mathfrak p_2, \dots, \mathfrak p_m$ are
  \marginpar{[213]} pairwise distinct prime ideals whose degrees are,
  respectively, $f_1,f_2,\dots,f_m$, and for each distinct prime
  function $P$ the corresponding prime ideal $\mathfrak p$ is the
  greatest common divisor of the ideals $\Of p$ and $\Of P(\theta)$.
\end{theorem}

\clearpage

\vspace{2\baselineskip}
\centerline{\S~3}
\vspace{2\baselineskip}

%\subsection{Paragraph 3: Determining when $p$ divides $k$}

From this theorem it follows that on the basis of a specific integer
$\theta$ from the field $\Omega$, which allows one to represent as
$\phi(\theta)$ infinitely many integers,\footnote{We are not sure what
  Dedekind means here; perhaps it is this. Given $\theta$, the
  infinitely many elements of $\Z[\theta]$ are algebraic integers in
  $\Omega$. The running assumption that $k\neq 0$ means that
  $\Z[\theta]$ has finite index in $\Of$.}  one can find the
factorization of all the prime numbers $p$ that do not not divide the
index corresponding to a the chosen $\theta$. It is therefore very
important to know whether a prime number $p$ is a divisor of the index
$k$ or not.\footnote{The goal of this section is to characterize the
  primes $p$ that divide the index of $\theta$. Dedekind points out
  that this is easy if we have the discriminant $D$ but his goal is to
  answer the question solely in terms of the minimal polynomial $F$.}
If we have a basis $\omega_1,\omega_2,\dots,\omega_n$ of the domain
$\Of $, or even just know the fundamental number $D$ of the field
$\Omega$, it is easy to answer the question, since in that case we can
find $k$ directly. From the coefficients of the equation $F(\theta)=0$
we can compute its discriminant
\[ \Delta(1,\theta,\theta^2,\dots,\theta^{n-1}) = (-1)^{\frac12n(n-1)}
  N(F'(\theta)) = Dk^2,\] and from that we can find the square of the
index $k$ by dividing by $D$.  In most investigations, however, things
are very different, since only the equation $F(\theta)=0$ is known,
and not the fundamental number $D$ of the corresponding field
$\Omega$. We would like to decide on that basis\footnote{That is,
  solely on the basis of the equation.} whether or not a specific
prime number $p$ divides the unknown index of the number $\theta$.
This is in fact possible, as we will now show, with the help of the
theory of higher congruences. Using our previous notation, the answer
turns out to depend on the nature of the function\footnote{This
  equation defines the polynomial $M$: it is the result of dividing
  the difference by $p$. Since $F$ is congruent mod $p$ to the
  product, $M$ is a polynomial with rational integer coefficients.}
$M$ that appears in the identity
\[ F = P_1^{e_1}P_2^{e_2}\dots P_m^{e_m} - pM.\]
This will be the content of the next two theorems.

\begin{theorem}[II] If the index of the number $\theta$ is not
  divisible by $p$, then $M$ cannot be divisible mod $p$ by any
  prime function $P$ whose square divides $F$ mod $p$.
\end{theorem}

To prove this,\footnote{\label{yikes}The theorem specifies a property
  of $M$ modulo $p$, but choosing different lifts for the factors
  $P_i$ can change $M$ (even mod $p$). So we need to check that the
  divisibility property we are looking for is independent of the
  chosen lifts.

  \textbf{Lemma:} Suppose $P,R,S,T\in\Z[t]$, $P\equiv R\pmod p$,
  $S\equiv T\pmod p$, $P$ is irreducible mod $p$, and that
  \[ F = P^e S - pM = R^eT - pN\]
  with $e\geq 2$. Then $M-N$ is divisible by $P$.

  \textbf{Proof of Lemma:} The equation
  \[ P^eS -pM=R^eT-pN\]
  gives
  \[ M-N = \frac1p(P^eS-R^eT),\]
  so we need to know $P^eS-R^eT \pmod{p^2}$. Writing $R=P+pX$, $T=S+pY$
  we get
  \[ P^eS-R^eT \equiv P^eS - (P^e + epP^{e-1}X)(S +pY)%\pmod{p^2}
    \equiv p(eP^{e-1}S + P^eY) \pmod{p^2}.\] Dividing by $p$ gives
  \[ M - N \equiv eP^{e-1}S+P^eY \pmod p,\]
  and since $e\geq 2$ we are done. \qed

  Dedekind does not prove this lemma; rather, he \emph{deduces} it
  from the fact that the question of whether the index is divisible by
  $p$ is independent of the choice of lifts. But he does say it can be
  checked directly, and it seems better to do that.} we can use the
results in the previous paragraph, which were all obtained under the
assumption that $p$ does not divide $k$.  \marginpar{[214]} Retaining
the same notation we used there,\footnote{So $S$ is the product of all
  the irreducible factors different from $P$; in particular, $S$ is
  not divisible by $P$.}  write $F\equiv SP^e\pmod p$, or
\[ F = SP^e-pM,\] and suppose $e\geq 2$. Then $p$ is divisible by
$\mathfrak p^2$, so $\mathfrak a$ is divisible by $\mathfrak p$ and
$\mathfrak b$ is not.\footnote{Since $\mathfrak a$ and $\mathfrak b$
  are relatively prime.} Therefore,\footnote{Since the principal ideal
  $\Of\rho$ is equal to $\fp\mathfrak b$, we know that
  $\rho=P(\theta)$ is divisible by $\fp$ only once, and that
  $\sigma=S(\theta)$ is not divisible by $\fp$.} $\mathfrak p^e$ is
the highest power of $\mathfrak p$ dividing
$S(\theta)P(\theta)^e = pM(\theta)$. Since $p$ is
divisible\footnote{Dedekind means $\fp^e$ is the highest power of
  $\fp$ that divides $p$, of course.} by $\mathfrak p^e$, it follows
that $M(\theta)$ cannot be divisible by $\mathfrak p$ and so
$M\not\equiv 0 \pmod{p,P}$, as claimed.\footnote{This concludes the
  proof of Theorem~II.}

It is also possible to prove the theorem without using the results in
the previous paragraph, in following indirect but equivalent
form:

\begin{theorem} If $F$ is divisible mod $p$ by the square of an
  irreducible polynomial $P$, so that $F=SP^e-pM$ with $e\geq 2$, and
  $M$ is divisible by $P$, then the index $k$ of the number $\theta$
  will be divisible by $p$.
\end{theorem}

Let the letters $\rho$, $\sigma$, $\eta$ have the same meanings as in
the previous paragraph, so that we set
\[ \rho = P(\theta), \qquad \sigma=S(\theta), \qquad
  \eta=\sigma\rho^{e-1}.\] Using the results of \S~1,\footnote{If $p$
  does not divide the index $k$, then two elements of $\Of'$ are
  congruent mod $p$ if and only if the corresponding polynomials are
  congruent mod $p,F$. In our case the polynomial will have degree
  less than $n$, so being congruent mod $p,F$ is equivalent to being
  congruent mod $p$.} the proof of our theorem will be complete if we
can show that the number $\eta=S(\theta)P(\theta)^{e-1}$ must be
divisible by $p$, since the function $SP^{e-1}$ is of degree lower
than $n$ and not $\equiv 0 \pmod p$.  To prove that $\eta$ is
divisible by $p$, it suffices to show that each power of a prime ideal
dividing $p$ also divides $\eta$ (D.~\S~163, B.~\S~25). To this end
set
\[\mu=M(\theta);\] 
consider the equation \[\sigma \rho^e=\eta\rho = p\mu.\] First, if
$\fp$ is a prime ideal dividing $p$ but not dividing $\rho$, then from
$\eta\rho=p\mu$ it follows at once that $\eta$ is divisible by the
highest power of $\fp$ dividing $p$.  Next, suppose $\fp$ divides both
$p$ and $\rho$. Since $S$ and $P$ are relatively prime
functions,\footnote{Dedekind doesn't say so, but he means relatively
  prime mod $p$.} there exist (C.~4) two functions $U$, $V$
\marginpar{[215]} such that the congruence
\[ SU + PV \equiv 1 \pmod p\] holds.
From that we get the numerical congruences\footnote{For the first one,
  we just plug in $\theta$; for the second, remember that $\fp$
  divides both $p$ and $\rho$.}
\[ \sigma U(\theta) + \rho V(\theta) \equiv 1 \pmod p\]
\[\sigma U(\theta)\equiv 1 \pmod \fp.\]
and it follows that $\sigma$ is not divisible by
$\fp$. Let\footnote{From here on we are basically computing $\fp$-adic
  valuations.}$\fp^h$, $\fp^r$, $\fp^m$ be the highest powers of $\fp$
dividing $p,\rho,\mu$, respectively. Since $\sigma\rho^e=p\mu$ and
$\eta=\sigma\rho^{e-1}$, we see that \[ er=h+m,\] and also that the
highest power of $\fp$ appearing in $\eta$ is equal
to \[(e-1)r=h+m-r.\] Since we want to show that $\eta$ is divisible by
$\fp^h$, it remains to prove that
\[m\geq r.\] Now we have to consider two cases.\footnote{The two cases
  are $r\geq h$ and $r\leq h$. Dedekind will use the assumption that
  $e\geq 2$ to handle the first case and the assumption that $M$ is
  divisible by $P$ to handle the second.} In the first, $r\geq h$, we
use the first assumption of our theorem, namely that $e\geq 2$. Then
$h+m=er\geq 2r$, and so $m-r\geq r-h\geq 0$, as claimed. In the second
case, $r\leq h$, we use the second assumption in our theorem, namely
that $M\equiv 0\pmod{p.P}$, i.e., $M\equiv PT \pmod p$. Therefore
$\mu\equiv \rho T(\theta)\pmod p$. Since $\rho$ is divisible by
$\fp^r$, it follows from this congruence that $\mu$ is also divisible
by $\fp^r$, and so that $m\geq r$, as we wanted to prove.

Now that we have proved Theorem~II in two different ways, we will also
show the correctness of the converse.

\begin{theorem}[III]
  If $M$ is not divisible mod $p$ by any prime function $P$ whose
  square divides $F$ mod $p$, the index $k$ of the number $\theta$ is
  not divisible by $p$.
\end{theorem}

The same theorem clearly can also be stated in the following form:

\begin{theorem} If the index $k$ of a number $\theta$ is divisible by
  $p$, there exists a prime function $P$ dividing $M$ whose square
  divides $F$ modulo $p$.
\end{theorem}

We present the proof\footnote{The structure of the proof is as
  follows. If $p|k$ then there exists a polynomial $\phi(t)\in\Z[t]$
  such that $\phi(t)\not\equiv 0 \pmod p$ but $\phi(\theta)$ is
  divisible by $p$ in $\Of$. We look at $A=\gcd(F,\phi)$ (over $\F_p$,
  but choose a monic lift of degree $<n$) and set $F=AB-pM$. Then we
  show that any prime divisor $P$ of $B$ in $\F_p[t]$ also divides
  $M$, therefore divides $F$, and we can then show it divides $A$ as
  well, so that $P^2|F$.  Factoring out the largest power of $P$ gives
  $F=P^eA'B'-pM$ with $P\nmid A'B'$, $e\geq 2$, $P|M$, which is what
  we want. Notice that the polynomial denoted by $M$ might change in
  the course of the argument.} of the latter form [of the theorem],
because the assumption that $k$ is divisible by $p$ is easier to use,
insofar as (according to \S~1) it implies the existence of a number
\[ \phi(\theta) = x_0 + x_1\theta + x_2\theta^2+\dots
  +x_{n-1}\theta^{n-1} \] which is divisible by $p$ but whose
coefficients $x_0,x_1,x_2,\dots,x_{n-1}$ are not all\footnote{We can
  say something a little stronger that will help below: if we had all
  but $x_0$ divisible by $p$, then $\phi(\theta)\equiv x_0\pmod p$ and
  so $x_0$ is also divisible by $p$. This means that $\phi(t)$ is not
  a constant mod $p$.} divisible by $p$.  Let us first denote by $A$
the greatest common divisor of $\phi(t)$ and $F$ modulo $p$. The
degree of $A$ is smaller than $n$, since $\phi$ has degree smaller
than $n$, and it is also not $\equiv 0\pmod p$. Write
\[ F= AB - pM,\] so that $B$ is not a
constant.\footnote{$\deg(B)=\deg(F)-\deg(A)=n-\deg(A)\neq 0$.} There
exist (C.~4) two functions $\phi_1$,$\phi_2$ such that
\[ \phi(t)\phi_1(t)+F(t)\phi_2(t)\equiv A(t) \pmod p.\] From this it
follows\footnote{Since $F(\theta)=0$ and $\phi(\theta)$ is divisible
  by $p$, it follows that $A(\theta)$ is divisible by $p$. In
  particular, $A$ cannot be a constant mod $p$. This transfers the
  assumption that $\phi(\theta)$ is divisible by $p$ to $A(\theta)$
  where $A|F$ in $\F_p[t]$.} that the number $A(\theta)$ is also
divisible by $p$,\footnoteB{In a similar way one can easily show that
  the criterion for the divisibility by $p$ of a number $\phi(\theta)$
  consists in the congruence $\phi(t)\equiv 0 \mod{p,K}$, where $K$ is
  a completely determined divisor of the function $F$ modulo $p$.}
From that\footnote{The asuumption that $\phi(\theta)$ is divisible by
  $p$ leads to the conclusion that $\frac1p A(\theta)$ is an algebraic
  integer, so it satisfies a monic equation with integer
  coefficients. Multiplying by a power of $p$ gives the equation
  below.} we get an equation of the form
\[ A(\theta)^s + ph_1A(\theta)^{s-1} + \dots + p^sh^s = 0,\] where
$h_1,h_2,\dots,h_s$ are rational integers (D.~\S~160; B.~\S~13). Since
the equation $F(\theta)=0$ is irreducible, this results in an equation
that holds identically\footnote{I.e., an equation in $\Z[t]$.} in the
variable $t$ of the form
\[ A^s + ph_1A^{s-1} + \dots + p^sh^s = FG,\]
which implies also the congruence
\[ A^s\equiv 0 \pmod{p,F}.\] Therefore\footnote{This is the key
  conclusion. Since $A^s$ is divisible by $F$ in $\F_p[t]$, every
  irreducible factor of $F$ must also divide $A$.} the function $A$
must be divisible modulo $p$ by every prime function that divides $F$
modulo $p$ (C.~5 and 6).  Now taking the equation above that is
satisfied by the number $A(\theta)$ and multiplying it by
$B(\theta)^s$, and recalling that $A(\theta)B(\theta)=p M(\theta)$ we
get\footnote{After dividing by $p^s$.}
\[ M(\theta)^s + h_1M(\theta)^{s-1}B(\theta) + \dots + h_sB(\theta)^s
= 0,\]
and therefore an identity of the form
\[ M^s + h_1M^{s-1}B + h_2M^{s-2}B^2 + \dots + h_sB^s = FH.\]
so\footnote{$B$ is a divisor of $F$ in $\F_p[t]$, so every term but
  the first is divisible by $B$ in $\F_p[t]$.}
$M^s \equiv 0 \pmod{p,B}$, which again implies that any prime function
dividing $B$ modulo $p$ must also divide $M$.  But we proved above
that $B$ is not a constant, so it has at least one prime divisor $P$,
which must then also divide $M$.\marginpar{[217]} Since $F$ is a
multiple of $B$ mod $p$, it must also divide $F$. But every prime
function dividing $F$ must divide $A$, as we showed above, so $P$ must
divide both $A$ and $B$, which shows $P^2$ must divide $F$, since
$F\equiv AB\pmod p$. So we have shown\footnote{Factoring out the
  highest powers of $P$ dividing $A$ and $B$ we get
  $F \equiv P^eA'B'\pmod p$, with $P\nmid A'B'$. By unique
  factorization in $\F_p[t]$, $A'B'$ is the rest of the factorization
  of $F$ and $M=\frac1p(F-P^eA'B')$ is (a possible choice for) the
  polynomial we are studying. We know that $e\geq 2$ and $P$ divides
  $M$. As we observed above, this property is independent of the
  choice of $M$, so it follows that no matter how we factor $F$ we
  will have $F=P^eS-pM$ with $P$ dividing $M$.}  that there is a prime
function contained in $M$ whose square is contained in $F$, which is
what we wanted to prove.

From II and III, the question of whether $p$ divides $k$ reduces to
looking at the factorization
\[ F = P_1^{e_1}P_2^{e_2}\dots P_m^{e_m} - pM\] of any function $F$
into prime functions modulo $p$.\footnote{All we need to check is
  there exists an $i$ such that $e_i\geq 2$ and $P_i$ divides $M$.}
In particular, if $F$ is not divisible by the square of any prime
function,\footnote{This is not the interesting case, since, as
  Dedekind's footnote points out, if $F$ is not divisible by the
  square of a polynomial mod $p$ it follows that $p$ does not divide
  the polynomial discriminant $\Delta=Dk^2$, so of course it does not
  divide $k$.}  so that all the exponents $e_1,e_2,\dots ,e_m$ are
equal to $1$,\footnoteB{This will be the case if and only if the
  discriminant $\Delta(1,\theta,\theta^2,\dots,\theta^{n-1})$ of the
  equation $F(\theta)=0$ is not divisible by $p$.} or when it
happens % \footnote{This is the interesting case.}
that none of the prime
functions whose squares divide $F$ are contained in $M$, then $k$ is
not divisible by $p$, and Theorem I from \S~2 applies.
% \footnote{Since
%   one can always determine the factorization of $p$ from the $p$-adic
%   factorization of a polynomial, it would be interesting to see if
%   Dedekind's result can be interpreted in terms of $p$-adic
%   factorizations.}
But if there is a prime function dividing $M$ whose square also
divides $F$, then $k$ is divisible by $p$ and the second proof of
Theorem~II shows that the factorization of the ideal $\Of p$ into
prime factors is \emph{different} from the one determined in
Theorem~I.\footnote{Dedekind doesn't explain why, but we think it
  might be this. If $F\equiv P^eS\pmod p$, $P$ does not divide $S$,
  and $\fp=\gcd(p,P(\theta))$, we expect that $\fp^e$ is a divisor of
  $p$. But in the second proof of Theorem~II we showed that if $p|k$
  the number $\eta = S(\theta)P(\theta)^{e-1}$ is divisible by $p$,
  but it is not divisible by $\fp^e$. Hence $\fp^e$ does not in fact
  divide $p$.}

To this result we add the following remark. If the functions $R_1$,
$R_2$, \dots, $R_m$ are congruent to the functions $P_1$, $P_2$,
\dots, $P_m$, then we have
\[ F = R_1^{e_1}R_2^{e_2}\dots R_m^{e_m} - pN\] and $N$ certainly does
not need to be congruent mod $p$ to $M$. On the other hand, the
divisibility of the index $k$ by $p$ is independent of the choice of
(lifts of) the divisors mod $p$, so we must have that the property of
$M$ that is key for this result will also hold for $N$.  This can
easily be confirmed directly by calculation.\footnote{We did this in
  the footnote on page~\pageref{yikes}. It is not quite clear that
  Dedekind's argument works without proving this first, but here he
  notes that it can easily be checked by a direct calculation.} If we
denote by $Q$ the product of all the prime functions contained in $F$
whose squares are \emph{not} contained in $F$, one can, by a suitable
choice of the functions $R_1, R_2,\dots, R_m$, always
arrive\footnote{If $F=PB\pmod p$ with $P\nmid B$ in $\F_p[t]$, we can
  always replace $P$ with $P+pC$ where $P\nmid C$ in $\F_p[t]$. That
  replaces $M$ by $M+CB$, which is not divisible by $P$ in $\F_p[t]$.}
at a function $N$ which is relatively prime to $Q$, but if there is a
prime function $P$ that divides $M$ \marginpar{[218]} such that $P^2$
divides $F$, a calculation shows that then $P$ divides $N$ as
well.\footnoteB{It follows from this that the ideal theory of
  Zolotareff is limited to the case in which the index $k$ is not
  divisible by $p$. At least this seems to follow from the following
  words, which we can find in the abstract mentioned above
  (\textit{Jahrbuch \'uber die Forstschritte der Mathematik}, Vol.~6):
  ``To present the theory in its simplest form, the author assumes
  that $F_1(x)$ is not divisible by any of the functions
  $V,V_1,V_2\dots$. If this condition does not hold, one can transform
  the equation $F(x)=0$ modulo $p$ so that it does hold. The author
  reserves the discussion of this transformation for another
  opportunity.'' --- Since according to my investigations (see \S~5 of
  this paper) there exist fields in which the indices of \emph{all}
  integers $\theta$ are divisible by a certain prime number $p$, it
  follows that \emph{all} equations $F(\theta)=0$ have the unfortunate
  property that impedes the application of Zolotareff's theory. Hence
  I suppose that there is a misunderstanding in the quoted words from
  the abstract. It is possible that the author's completion of the
  theory will be based on considerations similar to those in Selling's
  theory of ideal numbers (Schl\"omich's \textit{Zeitschirft},
  Vol.~10, p.~12ff.)}\footnote{The reference in Dedekind's footnote,
  in which ``12'' should read ``17,'' is \cite{Selling}.}

  % to the \textit{Zeitschrift f\"ur Mathematik and Physik}, then edited
  % by Schl\"omich. The article is \cite{Selling}.}

\vspace{2\baselineskip}
\centerline{\S~4}
\vspace{2\baselineskip}

%\subsection{Paragraph 4: A Criterion for Common Index Divisors}

In the number domains $\Of$ first considered by Kummer, which come
from a primitive root of the equation $\theta^m=1$, the happy
circumstance occurs that the powers
$1, \theta, \theta^2, \dots, \theta^{n-1}$, with $n=\phi(m)$, form a
basis\footnote{If $\theta$ is an $m$-th root of unity, we have
  $\Of=\Z[\theta]$, which is what allows Kummer's general approach to
  work.}  of the domain $\Of$. It follows that the index $k$ of the
number $\theta$, on which the entire investigation is based, is
$=1$. I soon realized, however, that in the general investigation of
\emph{any} finite field $\Omega$ and domain $\Of$ containing all the
integers in $\Omega$ this simple case rarely occurs.
\footnote{Dedekind says ``rarely'' or ``exceptionally.'' It is unclear
  how many examples he knew at this point, so this is an impressive
  insight. One expects that in fact the set of number fields with
  monogenic rings of integers has density zero. See, for example,
  \cite{monogenic}.}  I thought for a
\label{longtime} long time that it was likely that for every given
prime number $p$ it might be possible to find an integer $\theta$ in
the field $\Omega$ whose index is not divisible by $p$. If so, with
the help of that $\theta$ we could succeed in determining the ideal
factors of $p$. Since all my attempts to prove the existence of such a
number $\theta$ were unsuccessful, I finally decided, if possible, to
prove that this assumption was false. I achieved this goal, as I have
already indicated in the \emph{G\'ottingischen gelehrten Anziegen} of
September 20, 1871, through the considerations that form the content
of this and the following paragraph.\footnote{This is one of
  Dedekind's tantalizing accounts of why he ended up creating ideal
  theory: the ``local'' approach that might allow one to reduce
  everything to ``higher congruences'' is defeated by the existence of
  common index divisors. Kronecker \cite[\S~25]{Grundzuge} makes a
  similar argument.}

Let $p$ be a fixed prime number and let\footnote{Dedekind assumes he
  knows the factorization of $p$.} $\fp_1,\fp_2,\dots\fp_m$ be all the
distinct prime ideals dividing $p$; we will denote their degrees by
$f_1,f_2,\dots,f_m$, so that, for example, $N(\fp_1)=p^{f_1}$. If
there exists an integer $\theta$ whose index $k$ is not divisible by
$p$, it follows from Theorem~I that there exist $m$ polynomials
$P_1,P_2,\dots,P_m$ of degrees $f_1,f_2,\dots,f_m$, pairwise
\emph{incongruent} modulo $p$.\footnote{Since different prime ideals
  dividing $p$ correspond to distinct irreducible polynomials in
  $\F_p[t]$, such polynomials must exist. Conversely, if we cannot
  find enough irreducible polynomials of the required degrees, there
  cannot be any polynomial $F$ whose factorization matches the
  factorization of $p$, and so $k$ must be divisible by $p$ no matter
  which $\theta$ is chosen.} It is now of the greatest importance for
our investigation that this conclusion may be reversed,\footnote{So
  the existence of sufficiently many irreducible polynomials is enough
  to guarantee the existence of the appropriate $\theta$.}  so that
the following theorem holds.\footnote{On page 456 of \cite{Hasse},
  Hasse says that ``In deriving this criterion, Hensel gave the first
  demonstration of the power of his new foundation of algebraic number
  theory.'' The criterion there comes with an explicit formula for the
  number of irreducible monic polynomnials of degree $f$ in
  $\mathbb F_p[t]$, but this formula was certainly known to Dedekind
  as well. Hensel proves the same theorem in \cite{Hensel1894b}; see
  \cite{GW1}. Note, however, that to use this criterion one needs to
  know the factorization of $p$. That motivated Hensel to look for
  another criterion in \cite{Hensel1894b}, but the real solution came
  from the theory of $p$-adic numbers, which is what Hasse refers to
  as the ``new foundation of algebraic number theory.''}

\begin{theorem}[IV]
  Let $f_1,f_2,\dots,f_m$ be the degrees of the distinct prime ideals
  $\fp_1,\fp_2,\dots,\fp_m$ contained in $p$. Suppose that modulo $p$
  there exist $m$ incongruent prime functions $P_1,P_2,\dots,P_m$ of
  degrees $f_1,f_2.\dots, f_m$ respectively. Then there exists an
  integer $\theta$ in $\Omega$ whose index $k$ is not divisible by
  $p$.
\end{theorem}

Before giving a proof of this theorem, we will make a few general
observations\footnote{Making ``some observations'' is Dedekind's way
  to prove some lemmas. We will indicate each observation with a
  footnote.} that do not depend on all of its hypotheses.

Let $\fp$ be a prime ideal dividing $p$, of degree $f$ at least
one.\footnote{The first lemma says that if $\fp$ divides $p$ and has
  degree $f$, then there exists an irreducible polynomial
  $P\in\F_p[t]$ of degree $f$ and an element $\alpha\in\Of$ such that
  $P(\alpha)\equiv 0\pmod p$.}  Then \emph{all} the integers $\omega$
of the field $\Omega$ satisfy the congruence\footnote{Since $\Of/\fp$
  is a field with $p^f$ elements, every element has order $p^f-1$.}
\[ \omega^{p^f}-\omega\equiv 0\pmod{\fp}\]
(D.~\S~163;B.~\S~26, 3${}^{\text{o}}$).
Now if $t$ is a variable, the function
\[t^{p^f}-t\] is congruent mod $p$ to the product of all the
incongruent-mod-$p$ prime functions whose degree is a divisor of the
number $f$ (C.~19).\footnote{The (unique up to isomorphism) field with
  $p^f$ elements contains all the fields with $p^d$ elements with
  $d|f$. Again we notice that \cite{Abriss} is basically a theory of
  finite fields.} Among them we can choose \emph{at will} a prime
function $P$ whose degree if $=f$; this is always possible because
there always exists at least one such function (C.~20).\footnote{So
  $P$ is monic, irreducible mod $p$, and of degree $f$, and therefore
  a divisor of $t^{p^f}-t$.}  Then
\[ t^{p^f}-t \equiv P(t)H(t) \pmod p,\]
and therefore
\[ \omega^{p^f}-\omega \equiv P(\omega)H(\omega) \pmod p.\]
Since $\fp$ divides $p$, we see that
for \emph{every} number $\omega$ contained in $\Of$ we have the
congruence
\[ P(\omega)H(\omega)\equiv 0 \pmod{\fp}.\] Therefore the number of
roots $\omega$ that are incongruent modulo $\fp$ is exactly
$=[\Of:\fp]=N(\fp)=p^f$, therefore equal to the degree of the
congruence. Using the same simple arguments as in rational number
theory (D.~\S~26), one can easily prove that a congruence of degree
$r$ modulo a prime ideal $\fp$ can have no more incongruent roots in
the number domain $\Of$ than\marginpar{[220]} the degree $r$. I will
omit the proof for brevity.\footnote{A polynomial over a field cannot
  have more roots than its degree.}  Therefore in our case the congruence
$H(\omega)\equiv 0 \pmod\fp$ can have at most $(p^f-f)$ incongruent
roots, and it follows that the representatives $\omega$ of the $f$
other number classes must satisfy the congruence
$P(\omega)\equiv 0\pmod\fp$.  For our purposes, however, it is
sufficient to know that this congruence has at least one root. Let
$\alpha$ be one such root, so that
\[ P(\alpha)\equiv 0 \pmod\fp.\] We now consider all the number of the
form $\phi(\alpha)$ and we want to prove\footnote{The second lemma
  says that an element of $\Z[\alpha]$ is divisible by $\fp$ if and
  only if $\alpha=\phi(t)$ with $\phi(t)$ divisible by $P(t)$ in
  $\F_p[t]$. Essentially, $\alpha$ is a generator of the residue field
  $\Of/\fp$ and its minimal polynomial over $\F_p$ is $P$.}  that the
congruence
\[ \phi(\alpha)\equiv 0 \pmod{\fp}\]
is equivalent to the function congruence
\[\phi(t)\equiv 0 \pmod{p,P}.\]
Indeed, if the latter congruence holds, then also
\[ \phi(t)\equiv P(t)\psi(t)\pmod p,\]
and so
\[ \phi(\alpha)\equiv P(\alpha)\psi(\alpha) \pmod p,\]
and since both of the numbers $p$ and $P(\alpha)$ are divisible by
$\fp$, we get $\phi(\alpha)\equiv 0\mod \fp$. Conversely, if $\phi(t)$
is \emph{not} divisible by the prime funciton $P(t)$ then $\phi(t)$
and $P(t)$ will be relatively prime functions, and if follows that
there exist two functions $\phi_1(t),\phi_2(t)$ such that the
congruence
\[ \phi(t)\phi_1(t)+P(t)\phi_1(t)\equiv 1 \pmod p\]
holds (C.~5). Then we have
\[ \phi(\alpha)\phi_1(\alpha)+P(\alpha)\phi_1(\alpha)\equiv 1 \pmod
  p,\] and if follows that $\phi(\alpha)$ is \emph{not} $\equiv
0\pmod\fp$. So we have proved the claim above.

In the case\footnote{The third lemma is that if necessary we can
  change $\alpha$ to make sure $P(\alpha)$ is not divisible by
  $\fp^2$.} that $p$ is divisible by $\fp^2$, we also want to choose
the root $\alpha$ of the congruence $P(\alpha)\equiv 0\pmod\fp$ so
that the number $P(\alpha)$ is \emph{not} divisible by $\fp^2$.  This
is always possible:\marginpar{[221]} if $\alpha$ is a root of the
congruence $P(\alpha)\equiv 0 \pmod{\fp^2}$, then one can choose a
number $\lambda$ that is divisible by $\fp$ but not by $\fp^2$ and set
$\alpha'=\alpha+\lambda$. Then\footnote{The editors of \cite{Gesam}
  note here that here $P''(\alpha)$ really should be
  $\frac{P''(\alpha)}{2}$ and similarly for higher terms. They do not
  discuss whether those denominators will create trouble for the
  argument.} 
\begin{align*}
  P(\alpha')
  & =P(\alpha)+\lambda P'(\alpha)+\lambda^2P''(\alpha) +  \dots\\
  &\equiv \lambda P'(\alpha) \pmod{\fp^2}
\end{align*}
and since the derivative function $P'(t)$ has degree $(f-1)$ and is
not $\equiv 0 \pmod p$, it cannot\footnote{Dedekind knows that
  ``finite fields are perfect'' but only in the language of
  \cite{Abriss}.} be $\equiv 0\pmod{p,P}$, and therefore the number
$P'(\alpha)$ above is not divisible by $\fp$. So the number
$\lambda P'(\alpha)$, and therefore also the number $P(\alpha')$, is
divisible by $\fp$ but not be $\fp^2$. So we have proved the existence
of a number $\alpha'$. Let's remove the accent and thus assume that
$P(\alpha)$ is divisible by $\fp$ but not by $\fp^2$. \footnote{So now
  we have an irreducible polynomial $P$ of degree $f$, an element
  $\alpha\in\Of$ such that $P(\alpha)\equiv 0 \pmod{\fp}$; if $\fp^2$
  divides $p$, we can assume that $P(\alpha)\not\equiv0\pmod{\fp^2}$.}

Let\footnote{The fourth lemma is another translation into polynomial
  congruences, this time for divisibility by a power of $\fp$.}
$\fp^e$ be the highest power of $\fp$ contained in $p$; we want to
prove\footnote{The proof is identical to the previous one. The key
  result from \cite{Abriss} is that the gcd of two polynomials is a
  linear combination.} that the numerical congruence
\[ \phi(\alpha)\equiv 0\pmod{\fp^e}\]
is equivalent to the function congruence
\[\phi(t)\equiv 0 \pmod{p,P^e}.\]
If the latter holds, then
\[\phi(t) \equiv P(t)^e \psi(t)\pmod p,\]
so also
\[ \phi(\alpha) \equiv P(\alpha)^e\psi(\alpha) \pmod{p}.\] Since both
$p$ and $P(\alpha)^e$ are divisible by $\fp^e$, it follows that
$\phi(\alpha)\equiv 0 \pmod{\fp^e}$.  Conversely, if the function
congruence does \emph{not} hold, the greatest common divisor mod $p$
between $\phi(t)$ and $P(t)^e$ must be of the form $P(t)^s$ for some
$s<e$. So (C.~4) we have polynomials $\phi_1(t),\phi_2(t)$ such that
\[ \phi(t)\phi_1(t) + P(t)^e\phi_2(t) \equiv P(t)^s \pmod p.\]
Since both $p$ and $P(\alpha)^e$ are divisible by $\fp^e$, we get
\[ \phi(\alpha)\phi_1(\alpha)\equiv P(\alpha)^s \pmod{\fp^e}.\]
Since $s<e$ and $P(\alpha)$ is not divisible by $\fp^2$, it follows
that $\phi(\alpha)$ is \emph{not} $\equiv 0 \pmod{\fp^2}$, and our
claim is proved.

One can now apply the results described above\footnote{The final lemma
  shows, using the Chinese Remainder Theorem, that we can find a
  single $\theta$ satisfying the same properties as the
  $\alpha_i$. This will, of course, eventually be the $\theta$ whose
  existence is claimed in Theorem~IV.} to each of the prime
ideals $\fp_1,\fp_2,\dots,\fp_m$. One chooses \emph{arbitrary} prime
functions \marginpar{[222]} $P_1,P_2,\dots,P_m$ whose degrees
$f_1,f_2,\dots,f_m$ are the same as the degree of the corresponding
prime ideal. As above, one determines as many numbers
$\alpha_1,\alpha_2,\dots,\alpha_m$ such that
$P(\alpha_1),P(\alpha_2)\dots,P(\alpha_m)$ are respectively divisible
by $\fp_1,\fp_2,\dots,\fp_m$ and such that in the case that $p$ is
divisible by $\fp_r^2$ the corresponding $P_r(\alpha_r)$ is not
divisible by $\fp_r^2$.

Since the ideals $\fp_1,\fp_2,\dots,\fp_m$ are distinct, their squares
are [pairwise] relatively prime, so one can find (D.~\S~163;
B.~\S~26)\footnote{This is the Chinese Remainder Theorem in $\Of$;
  Dedekind does not use that name.} a number $\theta$ such that
\begin{align*}
  \theta&\equiv \alpha_1\pmod{\fp_1^2}\\
  \theta&\equiv \alpha_2\pmod{\fp_2^2}\\
  \dots&\dots\\
  \theta&\equiv \alpha_m\pmod{\fp_m^2}
\end{align*}
Then we have
\begin{align*}
  P_1(\theta)&\equiv P_1(\alpha_1)\pmod{\fp_1^2}\\
  P_2(\theta)&\equiv P_2(\alpha_2)\pmod{\fp_2^2}\\  
  \dots & \dots\\
  P_m(\theta)&\equiv P_m(\alpha_m)\pmod{\fp_m^2}\\
\end{align*}
It follows that the numbers
$P_1(\theta),P_2(\theta),\dots,P_m(\theta)$ are divisible respectively
by $\fp_1,\fp_2,\dots,\fp_m$, but, in the case when $p$ is divisible
by $\fp_r^2$, the number $P_r(\theta)$ is \emph{not} divisible by
$\fp_r^2$. The number $\theta$ therefore unites in itself all the
properties that each of the numbers $\alpha_r$ has with respect to the
corresponding prime ideal $\fp_r$.\footnote{So now we have a single
  number $\theta$ that ``works'' for all $i$ simultaneously, rather
  than individual $\alpha_i$.} Now\footnote{I think Dedekind has
  finished the ``general observations'' mentioned above, so that the
  proof of Theorem~IV now begins. One could, however, argue that this
  is one more lemma, and that the real proof begins when he invokes
  the key assumption in the next paragraph.} let
\[\Of p=\fp_1^{e_1}\fp_2^{e_2}\dots \fp_m^{e_m},\]
so that we have, by taking the norm,
\[ n=e_1f_1+e_2f_2+\dots +e_mf_m.\]
A number of the form $\phi(\theta)$ is divisible by one of the powers
$\fp_1^{e_1},\fp_2^{e_2},\dots,\fp_m^{e_m}$ if and only if the
corresponding function congruence
\begin{align*}
  \phi(t)&\equiv 0\pmod{p,P_1^{e_1}}\\
  \phi(t)&\equiv 0\pmod{p,P_2^{e_2}}\\
  \dots & \dots\\
  \phi(t)&\equiv 0\pmod{p,P_m^{e_m}}
\end{align*}
holds. An integer from the field is divisible by $p$
if and only if it is divisible by \emph{each} of the $m$ powers
$\fp_1^{e_1},\fp_2^{e_2},\dots,\fp_m^{e_m}$, and therefore a numerical
congruence
\[\phi(\theta)\equiv 0 \pmod p\]
is equivalent to the \emph{system} of $m$ function congruences
above.

So far \marginpar{[223]} we have intentionally put no restriction on
the \emph{choice} of the prime functions $P_1,P_2,\dots,F_m$ except
that their degrees are respectively those of the prime ideals
$\fp_1,\fp_2,\dots,\fp_m$, so that, for example, if $f_1=f_2$ nothing
stops us from choosing $P_1=P_2$.  We now want to introduce the main
assumption\footnote{The key assumption is used here.} of the theorem,
namely that we can find $m$ \emph{pairwise incongruent} prime
functions of the desired degree, and we will assume that
$P_1,P_2,\dots,P_m$ are such pairwise incongruent prime functions.
Then the powers $P_1^{e_1},P_2^{e_2},\dots,P_m^{e_m}$ will be pairwise
relatively prime; if we let
\[ R=P_1^{e_1}P_2^{e_2}\dots P_m^{e_m}\] be their product,
then\footnote{First use the equivalence we just proved.}  the
numerical congruence
\[ \phi(\theta) \equiv 0 \pmod p\] is equivalent the system of $m$
function congruences given above and so (C.~5)\footnote{This is the
  Chinese Remainder Theorem for polynomial congruences.} is equivalent
to the [single] function congruence
\[ \phi(t) \equiv 0 \pmod{p,R}.\]
Notice that the degree of the product $R$ is
\[e_1f_1+e_2f_2+\dots +e_mf_m,\] and so is
$=n$. Therefore\footnote{Since $\phi(t)$ is a polynomial of degree
  $n-1$ it can only be divisible by $R$ modulo $p$ if it is zero
  modulo $p$.}  a number
\[ \phi(\theta) = x_0 + x_1\theta + x_2\theta^2+\dots
  +x_{n-1}\theta^{n-1}\] can only be divisible by $p$ if
$\phi(t)\equiv 0 \pmod p$, i.e., only if all the $x_j$ are divisible
by $p$. It follows (from \S~1) that the index of $\theta$ is
\emph{not} divisible by $p$. So we have proved\footnote{If there are
  enough irreducible polynomials modulo $p$, we can find a number
  $\theta$ whose index is not divisible by $p$.} the Theorem stated
above, and we now want to add the following remark.\footnote{The
  remark is just that the polynomials we have found are exactly the
  factors (modulo $p$) of the irreducible polynomial of $\theta$.}

Since $k$ is not divisible by $p$, it is also not equal to $0$. Thus,
the number $\theta$ we obtained is the root of an irreducible
equation $F(\theta)=0$ of degree $n$. Then 
$F(\theta)\equiv 0 \pmod p$, so the function $F$ must be divisible by
$R$ mod $p$. Since both functions have degree $n$ and have highest
coefficient $1$, we must have $F\equiv R\pmod p$, i.e.
\[ F\equiv P_1^{e_1}P_2^{e_2}\dots P_m^{e_m} \pmod p,\] and we have
now returned to the starting point of our investigation in
\S~2,\footnote{To summarize, Dedekind has proved the following: if we
  know the factorization of $p$ in $\Omega$ and we choose \emph{any}
  list of polynomials of the appropriate degrees, then we can find an
  element $\theta$ that generates $\Omega$ over $\Q$ and whose index
  is not divisible by $p$. It will not necessarily be a nice
  generator.

  For example, suppose $\Omega=\Q(\sqrt2)$ and $p=7$. Then in fact
  $\Of=\Z[\sqrt2]$. The principal ideal $\Of(7)$ factors as
  $\Of 7=\Of(3+\sqrt2)\cdot\Of(3-\sqrt2)$. The factors
  $\fp_1=\Of(3+\sqrt2)$ and $\fp_2=\Of(3-\sqrt2)$ are both prime
  ideals of degree one. Let's deliberately make the ``wrong'' choice
  of two distinct polynomials of degree one in $\F_7[t]$: $P_1=t$,
  $P_2=t-1$. We now need to find $\alpha_1,\alpha_2$ such that
  $P_i(\alpha_i)\equiv 0 \pmod{\fp_i}$. Clearly we can take
  $\alpha_1=3+\sqrt2$ and $\alpha_2=4-\sqrt2$. Solving
  $\theta\equiv \alpha_i\pmod{\fp_i^2}$, we get
  $\theta\equiv25+27\sqrt2\pmod{49}$, so let $\theta=25+27\sqrt2$. The
  minimal polynomial for $\theta$ is
  \[ \theta^2-50\theta-833 \equiv \theta^2-\theta \pmod 7\] as
  expected, and the disciminant of $\theta$ is $5832$, so the index of
  $\Z[\theta]$ is $729=3^6$, which is not divisible by $7$.}

\vspace{2\baselineskip}
\centerline{\S~5}
\vspace{2\baselineskip}

% \subsection{Paragraph 5: A detailed example}
\marginpar{[224]} Our last investigation has yielded a criterion that
answers the question\footnote{More specifically, it gives an answer to
  the question of whether $\theta$ exists if we know the factorization
  of the principal ideal $\Of p$.} of whether $\Omega$ contains an
integer $\theta$ whose index is not divisible by $p$. When we
have \[\Of p=\fp_1^{e_1}\fp_2^{e_2}\dots \fp_m^{e_m},\] where
$\fp_1,\fp_2,\dots,\fp_m$ are distinct prime ideals whose degrees are,
respectively, $f_1,f_2,\dots,f_m$, then the singular case in which the
indices of \emph{all} the integers in $\Omega$ are divisible by $p$
happens when and only when it is not possible to find $m$ prime
functions of degree $f_1,f_2,\dots,f_m$ that are pairwise incongruent
mod $p$. Now we must ask whether the case when there are not enough
such polynomials ever does occur. To answer this, we will take
the simplest possible approach. The incongruent prime functions
\emph{of degree one} are the following:
\[ t, t+1, t+2, \dots, t+(p-1).\] Their number is $=p$. So the
singular case above will occur in a field $\Omega$ whenever a prime
number $p$ factors as the product of $p+1$ distinct prime ideals of
degree $1$.  By the norm computation above, the degree $n$ of such a
field must be\footnote{More generally, Dedekind's argument shows that
  if $p$ splits completely in a field of degree $n>p$, then it will be
  a common index divisor.} $=p+1$. If, in order to obtain the simplest
case, one takes the smallest prime number $p = 2$, the question arises
whether there are \emph{cubic} fields $\Omega$ in which the number $2$
is divisible by three distinct prime ideals of degree one.\footnote{So
  we can find fields with common index divisors if we can solve the
  following problem: given a prime $p$ and an integer $n>p$, find a
  field in which $p$ splits completely. Dedekind chooses $p=2$ and
  $n=3$.} In such a field, the indices of \emph{all} algebraic
integers will be \emph{even}. This investigation was carried out in
full generality the \textit{G\"ottingischen gelehrten Anzeigen} of
September 20, 1871, and led to an affirmative answer; here I will be
content to give a single example\footnote{Dedekind does much more than
  just give the example. Starting from an irreducible polynomial of
  degree three, Dedekind finds an integral basis, computes the
  discriminant, finds an explicit factorization of $2$, and computes
  explicitly the products of all the ideals that divide $2$. One gets
  the impression that he wants to show that everything in his theory
  can be computed explicitly, given enough time and parience.} that
has already been mentioned there.\footnote{The editors of \cite{Gesam}
  add a footnote: ``The notice [Anzeige] just mentioned contains an
  explanation of the method by which Dedekind came up with the example
  discussed here. It also contains another example of a field with a
  common index divisor, namely a quartic field in which the prime
  number $2$ decomposes into two prime ideals of degree two.''}

Let $\alpha$ be a root of the irreducible polynomial\marginpar{[225]}
of degree $3$\label{details}
\[ F(\alpha)=\alpha^3-\alpha^2-2\alpha-8=0.\] To find the
discriminant,\footnote{Dedekind uses the theorem that the discriminant
  of the minimal polynomial for $\alpha$ is equal to
  $(-1)^{n(n-1)/2} N(F'(\alpha))$.}  we find the number
\[ F'(\alpha) = \delta = -2 -2\alpha+3\alpha^2.\]
Then, repeatedly using $F(\alpha)=0$, we compute the products
\begin{align*}
  \delta\alpha &= 24 +4\alpha+\alpha^2\\
  \delta\alpha^2 &= 8+26\alpha+5\alpha^2
\end{align*}
and via linear elimination\footnote{We would say ``using the
  Cayley-Hamilton theorem.''}  of $1,\alpha,\alpha^2$ we find that
$\delta$ is a root of the equation
\[ \begin{vmatrix}-2-\delta & -2 & 3\\24 &4-\delta & 1\\8 &
  26&5-\delta \end{vmatrix}=0\]
that is,
\[ \delta^3 -7\delta^2-2012 = 0.\]
Therefore we have
\[ \Delta(1,\alpha,\alpha^2) = -N(\delta) = -2012 = -2^2\cdot 503.\]
Since $503$ is a prime number, the only two square divisors of this discriminant
are $1$ and $4$, so the index $k$ of the number $\alpha$ is either $1$
or $2$.
It is therefore the function
\[F(t)=t^3-t^2-2t-8\]
that we need to investigate modulo $p=2$. Clearly
\[ F = P_1^2P_2-2M\equiv P_1^2P_2 \pmod2,\]
where
\[ P_1=t,\quad P_2=t-1,\quad M=t+4.\]
Since $P_1$ is a factor of $M$ and $P_1^2$ is a factor of $F$ modulo
$2$, it follows\footnote{If we believe the theorems we are done: this
  is the condition for $2$ to be an index divisor for this polynomial,
  and so $k=2$. But Dedekind will check that his theorems are true
  each time.} (from the second proof of Theorem~II in \S~3) that the
number 
\[ P_1(\alpha)P_2(\alpha) = \alpha(\alpha+1)\]
is divisible by $2$, and therefore $k=2$.
That is immediately confirmed by the fact that the number
\[ \beta=\frac12\alpha(\alpha-1)-1\]
turns out to be an algebraic integer.
In fact, using the equation $F(\alpha)=0$ we find that
\begin{align*}
  \alpha^2&=2+\alpha+2\beta\\
  \beta^2&=-2+2\alpha-\beta\\
  \alpha\beta&=4
\end{align*}
and so
\[ \beta^3 + \beta^2+2\beta-8=0.\]
It \marginpar{[226]} follows that
\begin{align*}
  1&= 1\cdot 1 + 0\cdot\alpha+ 0\cdot\beta\\
  \alpha&=0\cdot 1 + 1\cdot\alpha + 0\cdot\beta\\
  \alpha^2&=2\cdot 1 + 1\cdot\alpha+2\cdot\beta
\end{align*}
and so
\[ \Delta(1,\alpha,\alpha^2)=\begin{vmatrix}1 & 0 & 0\\0 & 1 & 0\\2 &
1 &2 \end{vmatrix}^2\Delta(1,\alpha,\beta)=2^2\Delta(1,\alpha,\beta),\]
from which we see that
\[ \Delta(1,\alpha,\beta)=-503.\] Since this number is not divisible
by any square (except $1$), it is the fundamental number $D$ of our
cubic field $\Omega$, and the numbers $1,\alpha,\beta$ are a basis for
all the integers $\omega$ belonging to the domain $\Of$, i.e.,
\[ \Of=[1,\alpha,\beta]\]
in the notation we have used before.\footnote{We have computed the
  discriminant $D$ and an integral basis of $\Of$.}
Any algebraic integer in $\Of$ can be written in the form
\[ \omega = z + x\alpha + y\beta,\]
where $z,x,y$ are arbitrary rational integers.

We now want to use these results to determine the
factorization\footnote{This is the key: he wants to check that $2$
  splits completely, but he cannot use Theorem~I, so he goes for a
  direct computation. The first step is to show $2$ is unramified.} of
the number $2$. Since
\begin{align*}
  \alpha^2 &=2+\alpha+2\beta \equiv \alpha \pmod 2\\
  \beta^2 &= -2 +2\alpha + \beta \equiv \beta \pmod 2
\end{align*}
we get
\[(z+x\alpha+y\beta)^2 \equiv z^2 + x^2\alpha^2 + y^2\beta^2 \equiv
  z+x\alpha+y\beta \pmod 2,\] so that \emph{every} $\omega\in \Of$
satisfies $\omega^2-\omega\equiv 0 \pmod2$. If follows, first, that
$2$ cannot be divisible by the square of a prime ideal.\footnote{The
  next two sentences give a proof.} Indeed, if
$\Of(2)=\fp^2\mathfrak q$, with $\fp$ a prime ideal in $\Of$ or even
any ideal different from $\Of$, then $\mathfrak{pq}$ is not divisible
by $\Of(2)$, so there exists an element $\omega$ such that
$\mathfrak{pq}$ divides $\omega$ but $2$ does not divide
$\omega$. Then $\omega^2$ is divisible by $\fp^2\mathfrak{q}^2$ and
therefore by $2$, and this contradicts\footnote{Indeed,
  $0\equiv \omega^2\not\equiv \omega\pmod 2$.} the congruence
$\omega^2\equiv \omega\pmod2$ above.  So $\Of(2)$ is either a prime
ideal or a product of distinct prime ideals.\footnote{Now we need to
  prove $2$ is not a prime in $\Of$.} Let $\fp$ be a prime ideal
dividing $2$. Then we must have $\omega^2\equiv \omega\pmod{\fp}$ for
\emph{every} $\omega$ in $\Of$. The number of incongruent roots of
this congruence is $(\Of;\fp)=N(\fp)$, but the\marginpar{[227]} number
of roots cannot be larger than the degree of the congruence, so we get
$N(\fp)\leq 2$ and hence $N(\fp)=2$. So $\fp$ is a prime ideal and is
not all of $\Of$, since $N(\fp)>1$. Therefore every prime ideal
contained in $2$ is of degree one, and it follows, since $N(2)=2^3=8$,
that
\[ \Of(2)=\mathfrak{abc},\] where $\mathfrak a$, $\mathfrak b$,
$\mathfrak c$ are \emph{three} distinct prime ideals of degree
one.\footnote{We are done. But Dedekind will prove it again.} This
shows that the singular case described above does occur, and we will
check\footnote{Again, we already know this must happen, but Dedekind
  will check it explicitly.} that indeed the indices of \emph{all} the
number $\omega$ will be divisible by $2$. In fact, if we
set\footnote{Dedekind is just computing $\omega^2$.}
\begin{align*}
  z' &= z^2 + 2x^2-2y^2+8xy\\
  x' &= x^2 + 2y^2 +2xz\\
  y' &= 2x^2-y^2+2yz
\end{align*}
we have 
\[ \omega^2 = z' + x'\alpha+y'\beta,\] from which it follows that the
index of $\omega$ is equal to the determinant
\[ \begin{vmatrix} 1 & 0 & 0\\z & x &y\\z' &x'&y'\end{vmatrix} =
  xy'-x'y=2x^3-x^2y-xy^2-2y^3,
\]
which is always \emph{even}.\footnote{This will later be called
  computing the ``index form,'' especially in the Kronecker
  school. Hensel showed in \cite{Hensel1894b} that the index form has
  content $1$, i.e., it does not have an integer factor bigger than
  $1$. On the other hand, it may be that the \emph{values} of the
  index form are always divisible by some prime. Here the form is
  congruent mod $2$ to $(x^2-x)y-(y^2-y)x$ and hence is always
  $0\pmod 2$. In the same paper Hensel proved that $p$ is a common
  divisor of the values if and only if the form involves $u^p-u$, as
  here.}

In order to complete our example and to confirm the predictions
derived from general theory by calculation,\footnote{So we are going
  to check everything explicitly.} we want finally to represent the
ideals that appear here in the form of \emph{finite modules of rank
  three}\footnote{Dedekind say ``dreigliedrigen Moduln,'' which means
  something like ``triple modules'' or ``trinomial modules.''}
(D.~\S~161; B.\S~3), i.e., to determine these ideals by finding their
bases. These representations are as follows:\footnote{Dedekind doesn't
  explain how he computed these (it's easy enough), but he will check
  that these modules are indeed ideals and that their product is $2$.}
\begin{align*}
  \mathfrak a &= [2,\alpha, 1+\beta]\\
  \mathfrak b &= [2, 1+\alpha,\beta]\\
  \mathfrak c &= [2,\alpha,\beta].
\end{align*}
The system $\mathfrak a$ of all numbers of the form
\[ \alpha'=2z+\alpha z+(1+\beta)y,\]
where $z,x,y$ are arbitrary rational integers, indeed has the
fundamental properties of an ideal, namely:

I. The sums and differences of two numbers $\alpha'$ in the system
$\mathfrak a$ belong to the same system $\mathfrak a$.

II. Each product of a number $\alpha'$ from the system $\mathfrak a$
and a number $\omega$ from the domain $\Of$ is still a number from the
system $\mathfrak a$.

\marginpar{[228]}The first property is clear. To prove the second it
suffices to check that the product of each of the basis numbers $2$,
$\alpha$, $1+\beta$ of $\mathfrak a$ by each of the basis numbers $1$,
$\alpha$, $\beta$ of $\Of$ belongs to $\mathfrak a$. This is clear
right away for the five products
\[ 2.1, \alpha.1, (1+\beta).1, 2.\alpha,2.\beta=-2+2(1+\beta).\]
For the remaining four the same follows from the equations
\[ \alpha.\alpha=\alpha+2(1+\beta), \quad \alpha.\beta=2.2,\]
\[ (1+\beta)\alpha=2.2+\alpha,\quad (1+\beta)\beta=-2+2\alpha.\]
In the same way one can ckeck that the systems $\mathfrak b$ and
$\mathfrak c$ are ideals.

The \emph{Norm} $N(\mathfrak m)$ of an ideal $\mathfrak m$ is the
number $(\Of:\mathfrak m)$ of numbers that are incongruent mod
$\mathfrak m$ (D.~\S~163; B.~\S~20), which is equal to the determinant
of the expressions that give the basis numbers of $\mathfrak m$ as
linear combinations of the basis numbers of $\Of$ (D.~\S~161; B.~\S~4,
4${}^{\text{o}}$). So, for example,
\[ N(\mathfrak a)
  = \begin{vmatrix}2&0&0\\0&1&0\\1&0&1\end{vmatrix}=2,\] 
and in the same way
\[ N(\mathfrak b)=N(\mathfrak c)=2.\] When, however, the norm of an
ideal is a prime number, that ideal must necessarily be a prime
ideal,\footnote{This checks that the three factors are prime ideals.}
since in general we have
$N(\mathfrak a_1\mathfrak a_2)=N(\mathfrak a_1)N(\mathfrak
a_2)$. Therefore $\mathfrak{a,\, b,\, c}$ are prime ideals. Further,
they are pairwise distinct, since the number $\beta$ belongs to both
$\mathfrak b$ and $\mathfrak c$ but not to $\mathfrak a$, and the
number $\alpha$, which is contained in $\mathfrak c$, is not contained
in $\mathfrak b$. The number $2$ is contained in all three ideals and
so must also be contained in the product $\mathfrak{abc}$, so
$\Of(2)=\mathfrak{mabc}$, where $\mathfrak m$ is some ideal. But
computing the norm we get
\[ N(2)=8= N(\mathfrak{m}) N(\mathfrak{a}) N(\mathfrak{b})
  N(\mathfrak{c}) = 8N(\mathfrak{m}),\] therefore $N(\mathfrak m)=1$,
so $\mathfrak m=\Of$ and $\Of(2)=\mathfrak{abc}$.\footnote{We have
  proved the factorization a second time.} But we also want to check
this result, which follows from general theorems, by a direct
computation, i.e., through the actual \emph{multiplication} of the
ideals (D.~\S~165; B.~\S~12).\footnote{Dedekind's approach initially
  distinguished between two notions of divisibility. On the one hand,
  an ideal $\mathfrak a$ divides $\mathfrak b$ if
  $\mathfrak b\subset \mathfrak a$. On the other, $\mathfrak a$
  divides $\mathfrak b$ if there exists an ideal $\mathfrak c$ such
  that $\mathfrak b=\mathfrak{ac}$. In \cite{B} he says that
  proving these two notions are equivalent is ``the main difficulty of
  the theory,'' which he has overcome. He has established the
  factorization of $2$ using the first point of view. Now he will
  check it from the second point of view.}

By the \emph{product} $\mathfrak{ab}$ of two ideals we understand the
system of all products $\alpha'\beta'$ and all sums of such products
$\alpha'\beta'$, where $\alpha'$, $\beta'$ are any numbers belonging
respectively to the ideals $\mathfrak a$, $\mathfrak b$
\marginpar{[229]}(D.~\S163; B.~\S~22). Such a product is therefore
first a finite module whose basis numbers are all the products of each
basis nubmer of $\mathfrak a$ by each basis number of $\mathfrak
b$. In our case, then, $\mathfrak{ab}$ is the finite module whose
basis numbers are the nine products
\[ 2.2=4,\quad 2(1+\alpha)=2+2\alpha,\quad 2.\beta=2\beta,\]
\[\alpha.3=2\alpha,\alpha(1+\alpha)=2+2\alpha+2\beta,\quad
  \alpha\beta=4,\]
\[ (1+\beta).2=2+2\beta,\quad (1+\beta)(1+\alpha)=5+\alpha+\beta,\quad
  (1+\beta)\beta=-2+2\alpha.\] Of those nine number only three are
mutually \emph{independent} (D.~\S~159; B.~\S~4), so by the method I
have described in detail (B.~\S~4,6${}^{\text{o}}$), we can reduce
this module with nine generators\footnote{He says ``neungliedrigen
  Modul.''} to one with three
generators.\footnote{``Dreigliedrigen.''} Doing this very simple and
easy calculation one gets the six following equations:
\begin{align*}
  \mathfrak a^2 &= [4,\alpha,3+\beta]; &\mathfrak{bc}
  &=[2,2\alpha,\beta]\\
  \mathfrak b^2 &= [4,1+\alpha,\beta]; &\mathfrak{ca}
  &=[2,\alpha,2\beta]\\
  \mathfrak c^2 &= [4,2+\alpha,2+\beta]; &\mathfrak{ab}
  &=[2,2\alpha,1+\alpha+\beta]
\end{align*}
We now proceed in the same way. Multiplying each of those by
$\mathfrak{a,\,b,\,c}$ using the same method, we obtain the following
ten principal ideals:
\begin{align*}
\mathfrak{abc} &= [2,\,2\alpha,\,2\beta] = \Of(2)\\
\mathfrak{a^2c} &= [4,\,\alpha,\,2+2\beta] = \Of\alpha\\
\mathfrak{b^2c} &= [4,\,2+2\alpha,\,\beta] = \Of\beta\\
\mathfrak{ac^2} &= [4,\,2+\alpha,\,2\beta] = \Of(\alpha-2)\\
\mathfrak{bc^2} &= [4,\,2\alpha,\,2+\beta] = \Of(2-\beta)\\
\mathfrak{a^2b} &= [4,\,2\alpha,\,3+\alpha+\beta] = \Of(3+\alpha+\beta)\\
\mathfrak{ab^2} &= [4,\,2+2\alpha,\,1+\alpha+\beta] = \Of(1+\alpha+\beta)\\
\mathfrak{a^3} &= [8,\,4+\alpha,\,3+\beta] = \Of(3+2\alpha+\beta)\\
\mathfrak{b^3} &= [8,\,1+\alpha,\,4+\beta] = \Of(1+\alpha)\\
\mathfrak{c^3} &= [8,\,2+\alpha,\,2+\beta] = \Of(\alpha+\beta-4)
\end{align*}
The ten numbers $\mu$ to which these principal ideals
$\Of\mu=[\mu,\,\alpha\mu,\,\beta\mu]$ correspond are connected to each
other by the following easily checked relations:
\[
  \alpha(\alpha-2)(1+\alpha) =2^3;\qquad
  \alpha\beta=(\alpha-1)(1+\alpha+\beta_2^2\]
\[  (\alpha-2)(3+\alpha+\beta) = 2\alpha;\qquad
  \alpha(2-\beta)=2(\alpha-2)\]
\[  (\alpha-2)(3+2\alpha+\beta) =\alpha^2;\qquad
  \alpha(\alpha+\beta-4)=(\alpha-2)^2\]
\marginpar{[230]}This example, to which one add many
others, makes it clear that there are fields $\Omega$ in which the
indices of \emph{all} integers are divisible by the same prime number
$p$. This result is not a welcome one in some respects. Indeed, there
are many theorems in the theory of ideals that would be easy to prove
via the theory of higher congruences were in not for the fact that
Theorem~I in \S~2 requires the assumption that the index $k$ of the
integer $\theta$ not be divisible by $p$. We have now seen, however,
that in many cases this hypothesis cannot be satisfied no matter which
number $\theta$ we choose, and it follows that the approach suggested
by that Theorem will not work in full generality. For example, I
mention the following important theorem which I also used in the
\textit{G\"ottingischen gelehrten Anzeigen} of 20 Spetember 1971:

\textit{The fundamental number $D$ of a field $\Omega$ is the product
  of those and only those rational prime numbers $p$ that are
  divisible by the square of a prime ideal in that
  field.}\label{ramprimes} 

If there is an integer in $\Omega$ whose index is not divisible by the
prime number $p$, the truth of this result clearly follows very
easily\footnote{If $p$ divides $D$, then it divides $d(\theta)$ for
  any $\theta$, which means that the corresponding polynomial $F(x)$
  has a double root modulo $p$ and hence is divisible by the square of
  an irreducible polynomial. If we assume there is a $\theta$ whose
  index is not divisible by $p$, we can use Theorem~I, which tells us
  that $p$ is divisible by the square of a prime ideal. Conversely, if
  $p$ is divisible by the square of a prime ideal, $F(\theta)$ is
  divisible by the square of an irreducible polynomial mod $p$, and so
  $p|d(\theta)$. Since $d(\theta)=k^2D$ and $p\nmid k$ then $p|D$.}
from \S~2. But this obviously does not lead to a proof of the general
theorem, and it was only after several unsuccessful attempts that I
succeeded in finding a general proof. I must, however, reserve the
detailed development of this subject, in which the theorem itself will
be considerably generalized, for another occasion.

\subsection{Notes by \"Oystein Ore}

In Dedekind's collected works, this paper is followed by some two
pages of ``Erl\"auterungen zur vorstehenden Abhandlung'' signed by
\"Oysten Ore, one of the three editors. We present only some
highlights of what Ore has to say. Quotations are from
\cite[pp.~230--232]{Gesam}.

``The problem of generalizing the Kummer theory of ideals in
cyclotomic fields to general fields leads naturally to a definition of
ideals by means of higher congruences. Selling
(\textit{Zeitschr. f. Math n. Phys.}, vol 10, pp.~17--47 (1865))
already takes this path, and it is possible, using Galoisian
imaginaries and other auxiliary fields, to obtain a general theory of
ideals in Galois fields. The prime ideal decomposition of a prime
number $p$ is obtained from the factorization mod $p^a$ of the
defining equation in these auxiliary fields. A proof of the invariance
of these ideals, i.e., of their independence of the chosen defining
equation, is not clear.'' The Selling article is the same one cited by
Dedekind, \cite{Selling}. See Ore's comments below on factoring modulo
powers of $p$.

``As can be seen from the introduction,'' Dedekind had tried this
method as well,\footnote{This is even clearer from the discussion in
  the beginning of \S~4, where Dedekind says he thought for a long
  time that this would be possible.}  but then abandoned it in favor
of an abstract theory of ideals as presented in the second edition of
Dirichlet's \textit{Zahlentheorie}. This form of the theory does not
give us, however, an explicit way to determine the factorization of
given numbers in the field. Theorem~I solves that problem for primes
that do not divide the index, but the existence of common index
divisors (or common inessential discriminant divisors) blocks that
path in general.

The next few paragraphs focus on Zolotarev. In Zolotarev's first paper
\cite{Zolo1}, something like Theorem~I is used as the definition of
prime divisors, but of course this means it does not work for all
primes. Ore says that Zolotarev's second paper \cite{Zolo2} solved the
problem, but to do that had to abandon the approach based on higher
congruences. Ore explains Zolotarev's ``semi-local'' approach; there
are good expositions in \cite{Tcheb} and \cite{Engstrom}. Ore says he
will not get into all the alternative ways to lay the foundations.

Kronecker's approach based on forms gives a theoretically very simple
determination (``eine theoretisch besonders einfache Bestimmung'') of
the prime divisors of a rational prime.\footnote{Does Ore mean that
  Kronecker's theory yields an algorithm? Or that it is computable in
  theory but not in practice?} ``As was first shown in full generality
by Hensel\dots{} there is a complete analogue to Dedekind's theorem
for all prime numbers in this theory.'' The reference is to
\cite{Hensel1894a}, where Hensel shows that one can overcome the
existence of common index divisors by studying the
``Fundamentalgleichung.''\footnote{This is the polynomial in $t$ with
  coefficients in $\Z[x_0,x_1,\dots,x_{n-1}]$ that has the generic
  algebraic integer $\omega =x_0+x_1\omega+\dots+x_{n-1}\omega^{n-1}$
  and its conjugates as roots. Hensel proved that the discriminant of
  this equation is a homogeneous form in $n$ variables with content
  $D$, and also that the factorization modulo $p$ of the fundamental
  equation always gives the correct factorization of $p$ in the
  corresponding number field.}

``However, this solution to the problem does not provide any
information about the relationship between the properties of the field
equation and prime ideal decompositions, as is the case with
Dedekind's theorem. In Hensel's $p$-adic theory of algebraic numbers,
this gap is partially filled by showing that the decomposition of the
defining equation into irreducible $p$-adic factors corresponds to the
decomposition of $p$ into prime ideal powers. For the complete
determination of the prime ideal decomposition, one must also use
Kronecker's theory here. (See K.~Hensel: \textit{Theorie der
  algebraischen Zahlen I}, Leipzig 1908.)'' So even in 1930 Ore does
not see the $p$-adic solution as a complete answer. The reference is
to Hensel's \cite{HenselBook1}; there was never a volume II.

``It can be shown, however, that the difficulties of Dedekind's theory
can be completely eliminated if, instead of congruences (mod $p$) one
always considers congruences (mod $p^\alpha$) where $\alpha>\delta$ if
the discriminant of the corresponding equation is exactly divisible by
$p^\delta$.  The corresponding irreducible factors are then not
determined (mod $p^\alpha$), but rather (mod $p^{\alpha-\delta}$). The
common index divisors then completely lose their exceptional character
and one obtains a clear correspondence between prime ideal
decomposition and factors of the equation (O.~Ore, \textit{Math.\
  Ann.}, Vol.~96, pp.~315--352 (1926) and Vol.~97, pp.~569--598
(1927)). Furthermore, Dedekind's representation of the prime ideals in
the form $\beta=(p,\phi(\theta))$ can be recovered by a method that
shows great similarity to the determination of the series development
of algebraic functions (O.~Ore, \textit{Math.\ Ann.}, Vol.~99,
pp. 84--117 (1928)).'' It is surprising that Ore does not recognize
this as equivalent to Hensel's $p$-adic approach.

The theorem in section 4 gives a criterion for the existence of common
index divisors. Hensel gave a different one in \cite{Hensel1894b} (see
\cite{GW1} for a translation); Hensel's criterion is in terms of the
index form. In the same paper Hensel also proved Kronecker's
conjecture that if a number field $K$ has common index divisors then
there exists an extension field for which the values of the index form
do not have a common divisor.

Hensel's criterion implies that if $K$ is a number field of degree $n$
for which $p$ is a common index divisor then $p<\frac12 n(n-1)$. Ore
notes several improvements on this estimate: E.\ v.\ Zylinsky proved
\cite{Zylinsky}, using Dedekind's criterion, that in fact
$p<n$. M.~Bauer showed \cite{Bauer} that if $p<n$ there always exists
a field of degree $n$ for which $p$ is a common index divisor. Bauer's
result also follows from general theorems of Hasse
\cite{HassePrescribed} showing the existence of fields in which $p$
has prescribed factorization.

\section{Conclusion}

Richard Dedekind was the first to give an example of a field in which
there is a common index divisor, in his 1871
\textit{Anzeige}. (Kronecker says he knew an example in 1858, but he
did not mention it in print until 1882.) In his 1878 paper, Dedekind
showed that such common index divisors were entirely a ``small
prines'' effect. Specifically, $p$ is a common index divisor for $K$ if and
only if translating its factorization in $\Of_K$ into a polynomial
factorization modulo $p$ requires too many irreducible
polynomials. This is exactly the criterion that Hasse gives in
\cite[p.~456]{Hasse} and attributes to Hensel.

As Ore pointed out in his notes, the usefulness of the criterion is
limited in that it requires knowledge of the factorization of
$p$. This is frustrating, since starting point of the investigation
was exactly the use of Dedekind's theorem to determine the
factorization of $p$. The status of the problem in 1878 was this:
\begin{itemize}
\item Given a generator $\alpha$ and a prime $p$, one can tell,
  looking at the factorization of the minimal polynomial modulo $p$,
  whether $p$ divides the index $(\Of:\Z[\alpha])$. This is one of the
  main results in the paper translated above.
\item If $p$ did not divide the index, one could determine its
  factorization in $\Of$ from the factorization modulo $p$ of the
  irreducible polynomial.
\item In some cases, however, it is impossible to find such a
  generator $\alpha$. This happens when there are not enough
  irreducible polynomials modulo $p$ to reflect the correct
  factorization of $p$ in $\Of$.
\item For such common index divisors, no algorithm was available to
  determine the factorization.
\end{itemize}

Kronecker assigned the problem to Hensel sometime in the early
1880s. It was the topic of Hensel's dissertation and of several papers
until the culminating paper of 1894, which we translate in
\cite{GW1}. This suggests that neither Kronecker nor Hensel had fully
absorbed Dedekind's 1878 paper, though both cite it.

%\bibliographystyle{plain}
%\bibliography{hensel}

\end{document}